\documentclass[11pt]{article}

\usepackage[skip=0.4\baselineskip]{parskip}

\usepackage{amsfonts,amssymb,amsthm,amsmath}
\usepackage{mathtools}
\usepackage{stmaryrd}

\usepackage{xcolor}
\usepackage{tikz}
\usepackage{tikzstyle}
\usetikzlibrary{positioning, arrows.meta,calc}
\usetikzlibrary{decorations.pathreplacing,decorations.markings}
\usetikzlibrary{matrix}
\usetikzlibrary{decorations}

\usepackage{graphicx}
\usepackage{fullpage}
\usepackage[shortlabels]{enumitem}

\numberwithin{equation}{section}

\usepackage{url}
\usepackage[linktoc=all,colorlinks,unicode=true]{hyperref}
\usepackage[capitalise,compress,nameinlink,noabbrev]{cleveref}

\hypersetup{
  linkcolor={blue!70!black},
  citecolor={blue!70!black},
  urlcolor={blue!70!black}
}

\usepackage[
  backend=biber,
  style=numeric,
  sortcites=true,
 giveninits=true,
  url=false,
  eprint=true,
  maxbibnames=99,
  doi=true
]{biblatex}
\usepackage{thmtools}
\newtheorem{thm}[equation]{Theorem}
\newtheorem{cor}[equation]{Corollary}
\newtheorem{lem}[equation]{Lemma}
\newtheorem{prop}[equation]{Proposition}
\newtheorem{claim}[equation]{Claim}

\newtheorem{ques}[equation]{Question}
\newtheorem{prob}[equation]{Problem}

\theoremstyle{definition}
\newtheorem{defn}[equation]{Definition}

\newtheorem*{funding}{Funding}
\newtheorem*{ai}{AI Declaration}

\theoremstyle{remark}

\newtheorem{rem}[equation]{Remark}

\crefname{thm}{Theorem}{Theorems}
\Crefname{thm}{Theorem}{Theorems}
\crefname{lem}{Lemma}{Lemmas}
\Crefname{lem}{Lemma}{Lemmas}

\newenvironment{poc}{\begin{proof}}{\end{proof}}

\title{Sidorenko Inequalities in Oriented Graphs}
\author{Veronica Bitonti\thanks{Mathematical Institute, University of Oxford. E-mail: {\tt veronica.bitonti@maths.ox.ac.uk}, {\tt emma.hogan@maths.ox.ac.uk}, {\tt dmitry.tsarev@maths.ox.ac.uk}.}\protect\footnotemark[1]
\and Emma Hogan\protect\footnotemark[1] \and Jonathan A. Noel\thanks{Department of Mathematics and Statistics, University of Victoria. E-mail: {\tt noelj@uvic.ca}.} \and Dmitry Tsarev\footnotemark[1]}

\DeclareTextCompositeCommand{\v}{OT1}{l}{l\nobreak\hspace{-.1em}'}
\DeclareTextCompositeCommand{\v}{OT1}{t}{t\nobreak\hspace{-.1em}'\nobreak\hspace{-.15em}}

\DeclareMathOperator{\rng}{rng}

\begin{document}

\maketitle

\begin{abstract}
We introduce relative Sidorenko and relative anti-Sidorenko inequalities for oriented graphs, strengthening the corresponding well-studied notions for tournaments by allowing arbitrary host graphs. An orientation $\vec H$ of a graph $H$ is \emph{relatively Sidorenko} if
$\hom(\vec H,\vec G) \ge 2^{-e(H)}\hom(H,G)$
for every oriented graph $\vec G$ with underlying graph $G$, and \emph{relatively anti-Sidorenko} if the reverse inequality holds. We provide a complete characterisation of the class of connected relatively Sidorenko graphs in terms of weakly norming graphs, and exhibit several families of relatively Sidorenko and relatively anti-Sidorenko oriented graphs. Despite the substantially stronger definitions, many known tournament Sidorenko and anti-Sidorenko examples satisfy the relative inequalities, suggesting a broader theory of these phenomena in arbitrary host graphs and motivating several new structural and classification problems.
\end{abstract}

\section{Introduction}

A major theme in extremal combinatorics is identifying substructures whose density in a large combinatorial object is minimised or maximised by a random construction. A quintessential example is Sidorenko's conjecture~\cite{Sidorenko93}, which informally states that, for any bipartite graph $H$ and $p \in (0,1)$, the number of labelled copies of $H$ in a large graph $G$ of edge-density $p$ is asymptotically minimised by taking $G$ to be a random graph of density $p$. Following decades of work on Sidorenko’s conjecture, OpenAI recently claimed to have disproved it, exhibiting a counterexample on $35$ vertices and $66$ edges \cite{OpenAISidorenko2026}. Given this, it is interesting to better understand the barrier between bipartite graphs and other structures that satisfy Sidorenko-type inequalities and those that do not. 

Sidorenko's conjecture belongs to a broader class of homomorphism-counting problems in extremal graph theory, many of which are closely connected to the study of quasirandomness. One natural question is whether analogous extremal phenomena occur for oriented graphs. We first introduce some necessary terminology. An \emph{oriented graph} is a digraph with no loops and at most one arc between each pair of distinct vertices. We will always denote an oriented graph with an arrow, as in $\vec{H}$, and an undirected graph without an arrow. A \emph{homomorphism} from an oriented graph $\vec{H}$ to an oriented graph $\vec{G}$ is a map $\varphi:V(\vec{H})\to V(\vec{G})$ such that if $(u,v)$ is an arc of $\vec{H}$, then $(\varphi(u),\varphi(v))$ is an arc of $\vec{G}$. Let $\hom(\vec{H},\vec{G})$ denote the number of homomorphisms from $\vec{H}$ to $\vec{G}$. Homomorphisms and the function $\hom(H,G)$ for undirected graphs are defined analogously. 

Inspired by Sidorenko's conjecture, a natural question is whether there exist oriented graphs $\vec{H}$ such that, for a given host graph $G$, the number of homomorphisms from $\vec{H}$ into an orientation of $G$ is asymptotically minimised or maximised by a uniformly random orientation of $G$. This concept was first formally introduced when $G$ is a complete graph $K_n$, i.e. when $\vec{G}$ is a \emph{tournament}, by Zhao and Zhou~\cite{ZhaoZhou20}; we use the terminology from~\cite{FoxHimwichManiZhou25}.

\begin{defn}
\label{def:tournamentSid}
    An oriented graph $\vec{H}$ is \emph{tournament Sidorenko} if 
    \[\hom(\vec{H},\vec{K}_n) \geq (1/2)^{e(H)}(1 - o(1))n^{v(H)}\]
    for every tournament $\vec{K}_n$. Similarly, $\vec{H}$ is \emph{tournament anti-Sidorenko} if 
    \[\hom(\vec{H},\vec{K}_n) \leq (1/2)^{e(H)}n^{v(H)}\]
    for every tournament $\vec{K}_n$.
\end{defn}

This is asymptotically equivalent to saying $\vec{H}$ is tournament Sidorenko (resp. tournament anti-Sidorenko) if the number of labelled copies of $\vec{H}$ in a tournament $\vec{K}_n$ is minimised (resp. maximised) by the uniformly random orientation of $K_n$. The $o(1)$ term is present to account for the fact that $K_n$ does not have loops, and so there exist maps from $V(H)$ to $V(K_n)$ that are not graph homomorphisms, namely, the maps that send each endpoint of some edge in $H$ to the same vertex of $K_n$. It is straightforward to check that $\hom(H,K_n) = (1 - o(1))n^{v(H)}$. These notions are closely related to tournament quasirandomness. Roughly speaking, an oriented graph $\vec{H}$ \emph{forces quasirandomness} in tournaments if it satisfies one of the bounds in \cref{def:tournamentSid} and, whenever the bound is close to being tight, it implies that the global structure of $\vec{K}_n$ is similar to that of a uniformly random tournament.

There is now extensive literature on tournament Sidorenko and tournament anti-Sidorenko oriented graphs~\cite{ZhaoZhou20,SahSawhneyZhao23,FoxHimwichManiZhou24+,HeManiNieTungWei25+,ChenClemenNoel26+,ChenClemenNoelSharfenberg26+}, as well as tournament quasirandomness~\cite{ChungGraham91b,KalyanasundaramShapira13,CoreglianoRazborov17,CoreglianoParenteSato19,BucicLongShapiraSudakov21,GrzesikKralLovaszVolec23,Hancock+23,NoelRanganathanSimbaqueba26,Kral+26+}.
However, little is known about similar inequalities when the host graph is not 
required to be complete. In this paper, we remove the completeness assumption on the host graph and ask for inequalities that hold uniformly over all orientations of an arbitrary graph $G$. Formally, we define the following. 

\begin{defn}
\label{def:main}
An oriented graph $\vec{H}$ is \emph{relatively Sidorenko} if
    \[\hom(\vec{H},\vec{G})\geq (1/2)^{e(H)}\hom(H,G)\]
for every oriented graph $\vec{G}$, where $G$ is the undirected graph underlying $\vec{G}$. Similarly, $\vec{H}$ is \emph{relatively anti-Sidorenko} if
    \[\hom(\vec{H},\vec{G})\leq (1/2)^{e(H)}\hom(H,G)\]
for every oriented graph $\vec{G}$.
\end{defn}

Although this terminology is new, the underlying concepts are not. They appear implicitly in the work of Griffiths~\cite{Griffiths13} on quasirandom oriented graphs. Additionally, the property of being relatively Sidorenko is mentioned in~\cite[Section~5]{FoxHimwichManiZhou25}, where it is referred to as the ``second directed Sidorenko property''. In particular, Griffiths established the following examples.

\begin{thm}[Griffiths~{\cite[Proposition~3.1]{Griffiths13}}]
\label{th:GriffithsSid}
The oriented $4$-cycle in which every vertex is either a sink or a source is relatively Sidorenko, and the oriented $4$-cycle in which three edges are oriented clockwise and one edge is oriented anticlockwise is relatively anti-Sidorenko.
\end{thm}

Our main result is a characterisation of relatively Sidorenko oriented graphs in terms of the well-studied class of weakly norming graphs. The precise definition of weakly norming graphs is somewhat technical and is postponed to \cref{sec:Sid}, although we note that weakly norming graphs are necessarily bipartite \cite{Hatami10}. We say that an orientation of a bipartite graph $H$ is \emph{one-way} if there is a bipartition $(X,Y)$ of $H$ such that every edge of $H$ is oriented from $X$ to $Y$, and that a bipartite graph $H$ is \emph{reversible} if there is a bipartition $(X,Y)$ of $H$ and an automorphism $\phi$ of $H$ such that $\phi(X)=Y$. In other words, a connected graph $H$ is reversible if it has a unique one-way orientation, up to isomorphism. We prove the following characterisation result in \cref{sec:Sid}.

\begin{restatable}{thm}{relSidCharacterisation}\label{th:relSidCharacterisation}
Let $\vec{H}$ be a connected oriented graph with at least one arc. Then $\vec{H}$ is relatively Sidorenko if and only if it is a one-way orientation of a reversible weakly norming graph. 
\end{restatable}

While \cref{th:relSidCharacterisation} relates relatively Sidorenko oriented graphs to weakly norming graphs, we note that obtaining a structural characterisation of weakly norming graphs themselves remains an open problem. However, we derive several explicit families of relatively Sidorenko oriented graphs from \cref{th:relSidCharacterisation} and known results on weakly norming graphs; additional examples may be derived from the results of~\cite{ConlonLee17}.

\begin{restatable}{thm}{relSidExamples}\label{th:relSidExamples}
The one-way orientations of the following classes of graphs are relatively Sidorenko:
\begin{enumerate}[(a)]
    \item\label{cycles} the even cycle $C_{2k}$ for $k \geq 2$,
    \item\label{Knn} the complete bipartite graph $K_{k,k}$ for $k \geq 1$,
    \item\label{Knn-M} the complete bipartite graph $K_{k,k}$ with a perfect matching deleted for $k \geq 1$, and
    \item\label{Qd} the $k$-dimensional hypercube $Q_k$ for $k \geq 1$.
\end{enumerate}
\end{restatable}

Next, in~\cref{sec:antiSid}, we establish several classes of relatively anti-Sidorenko oriented graphs. The following theorem follows from an iterative construction for families of relatively anti-Sidorenko graphs, and demonstrates that the family of relatively anti-Sidorenko graphs is vastly different from the family of relatively Sidorenko graphs.
\begin{restatable}{thm}{antiSidChrom}\label{th:antiSidChrom}
For every $k\geq1$ there exists a relatively anti-Sidorenko oriented graph whose underlying graph has chromatic number $k$.
\end{restatable}
We additionally apply entropy methods to prove relative anti-Sidorenko inequalities for a certain class of highly symmetric trees that we call \emph{perfectly balanced}. A precise definition is given in \cref{sec:antiSid}, along with examples in \cref{fig:balancedTrees}.

\begin{restatable}{thm}{perfectlyBalanced}
\label{th:antiSidTrees}
Every perfectly balanced oriented tree is relatively anti-Sidorenko.
\end{restatable}

Fox, Himwich, Mani and Zhou~\cite{FoxHimwichManiZhou24+} conjectured that every undirected tree has at least one tournament anti-Sidorenko orientation. This conjecture remains open; see~\cite{HeManiNieTungWei25+,ChenClemenNoel26+,ChenClemenNoelSharfenberg26+} for recent progress. In contrast, in~\cref{sec:notAntiSid}, we show that not all trees admit relatively anti-Sidorenko orientations. In particular, we consider paths, all of which were shown  by Sah, Sawhney and Zhao to admit a tournament anti-Sidorenko orientation~\cite{SahSawhneyZhao23}. The following result shows that the situation in general host graphs is markedly different.

\begin{restatable}{thm}{notAntiSidPath}\label{th:path}
    No path with at least $3$ edges admits a relatively anti-Sidorenko orientation.
\end{restatable}

Furthermore, Fox, Himwich, Mani and Zhou~\cite[Theorem~1.11]{FoxHimwichManiZhou24+} prove that an oriented star is tournament anti-Sidorenko if and only if the in-degree and out-degree of the central vertex differ by at most one. Once again, the following result shows that the class of relatively anti-Sidorenko oriented stars is more constrained.

\begin{restatable}{prop}{antiSidStar}\label{prop:star}
    An oriented star with at least three vertices is relatively anti-Sidorenko if and only if the central vertex has in-degree equal to its out-degree.
\end{restatable}

The results in this paper leave open several classification questions for relatively Sidorenko and relatively anti-Sidorenko orientations. We conclude in~\cref{sec:concl} with some open problems.

\textbf{Related Independent Work:}
The notions of relatively Sidorenko and relatively anti-Sidorenko oriented graphs were independently studied by Nie and Wei~\cite{NieWei26+} under the names \emph{oriented Sidorenko} and \emph{oriented anti-Sidorenko}. In particular:
\begin{itemize}
\item The characterisation of connected relatively Sidorenko oriented graphs in our~\cref{th:relSidCharacterisation} is the same as~\cite[Theorem~1.2]{NieWei26+}. The proofs are similar, but differ in how reversibility of the underlying graph is established.
\item The relatively Sidorenko families in our~\cref{th:relSidExamples} also appear in the discussion following~\cite[Theorem~1.2]{NieWei26+}.
\item Our classification of relatively anti-Sidorenko stars in~\cref{prop:star} and our result on paths in~\cref{th:path} correspond to~\cite[Theorems~1.9 and~1.10]{NieWei26+}, respectively.
\end{itemize}

There are also substantial differences between our paper and~\cite{NieWei26+}.  We prove that every perfectly balanced oriented tree is relatively anti-Sidorenko (\cref{th:antiSidTrees}), together with a stronger geometric-mean inequality over all orientations of its underlying tree (\cref{th:antiSidTreesStronger}). We also construct relatively anti-Sidorenko oriented graphs with underlying chromatic number exactly $k$ for every $k\geq1$ (\cref{th:antiSidChrom}). The non-bipartite example constructed in the proof of~\cite[Proposition~1.8]{NieWei26+} is obtained from the second construction in~\cref{rem:antiSidConstruction} by taking $\vec{F}$ to be a single isolated vertex.

Nie and Wei additionally classify relatively anti-Sidorenko orientations of cycles, showing that Griffiths's orientation of $C_4$ is the only example~\cite[Theorem~1.11]{NieWei26+}. They further investigate when asymptotic equality in these inequalities forces quasirandom edge directions, a property they call \emph{symmetry forcing}~\cite[Theorem~1.3, Theorem~1.7 and Proposition~1.12]{NieWei26+}. They also establish results on undirected step forcing and strict convexity of homomorphism density functionals on non-negative kernels~\cite[Theorem~1.4, Corollary~1.5 and Theorem~1.6]{NieWei26+}.

\textbf{Notation and terminology:} 
Throughout the paper we use the following notation. For integers $n\geq1$, we let $[n]\coloneq \{1,\dots,n\}$. Given an undirected graph $H$, we denote by $V(H)$ and $E(H)$ the vertex set and edge set of $H$, respectively, and write $e(H) \coloneq |E(H)|$ and ${v(H) \coloneq |V(H)|}$. For an oriented graph $\vec{H}$, we denote by $A(\vec{H})$ the arc set of $\vec{H}$, and write $a(\vec{H}) \coloneq |A(\vec{H})|$. Given an oriented graph $\vec{H}$, we denote by $H$ the graph obtained by replacing each arc of $\vec{H}$ with an undirected edge, and we denote by $\vec{H}^{\text{rev}}$ the graph obtained by switching the direction on every arc of $\vec{H}$.

For a vertex $v$ in an undirected graph $H$, let $d_H(v)$ denote the \emph{degree} of $v$ in $H$; that is, the number of edges incident to $v$.  The \emph{in-degree} $d^-_{\vec{G}}(v)$ of a vertex $v$ in an oriented graph $\vec{G}$ is the number of arcs $(u,v)$ for some vertex $u \in V(\vec{G})$, and the \emph{out-degree} $d^+_{\vec{G}}(v)$ of $v$ is the number of arcs $(v,u)$ for some vertex $u \in V(\vec{G})$. A vertex with in-degree zero is called a \emph{source}, and a vertex with out-degree zero is called a \emph{sink}; an isolated vertex is both a source and a sink. Given a subset $S\subseteq V(\vec{H})$ of vertices of some oriented graph $\vec{H}$, the subgraph of $\vec{H}$ \emph{induced} by $S$ is the oriented graph $\vec{H}[S]$ with vertex set $S$ and all arcs $(u,v) \in A(\vec{H})$ such that $u,v\in S$. We similarly define the induced subgraph $H[S]$ of an undirected graph $H$.

\subsection*{Declarations}

\begin{funding}
The first author was supported by EPSRC grant EP/Z534870/1, the second author was supported by EPSRC grant EP/X013642/1, the third author was supported by NSERC Discovery Grant RGPIN-2021-02460 and the fourth author was supported by EPSRC grant EP/W523781/1.
\end{funding}

\begin{ai}
We proved \cref{th:relSidExamples,th:antiSidTrees,th:path}, \cref{prop:star} and various necessary conditions for relatively Sidorenko oriented graphs which are implied by \cref{th:relSidCharacterisation} without AI assistance. Our original proof of \cref{th:relSidExamples} used ideas from the literature on weakly norming graphs, particularly cut involutions and Cauchy--Schwarz trees from~\cite{ConlonLee17}, and motivated us to ask whether every one-way orientation of a connected reversible weakly norming graph is relatively Sidorenko in the conclusion of a draft of our paper. We provided this draft to ChatGPT 5.5, prompting it to investigate this question; it answered in the affirmative and supplied a proof. We then guided it, through a series of prompts, toward proving the converse, and it was successful in doing so. These AI-generated arguments informed the initial proof of \cref{th:relSidCharacterisation}, which has now been modified and completely rewritten by the authors for readability. The results on anti-Sidorenko oriented graphs were obtained without AI assistance. ChatGPT also assisted with copy-editing, produced TikZ figures, and suggested references to be checked by the authors. The authors take full responsibility for all content in the paper, including the accuracy of the proofs and the exposition.
\end{ai}

\section{Relatively Sidorenko Oriented Graphs}
\label{sec:Sid}

In this section we prove \cref{th:relSidCharacterisation,th:relSidExamples}. Before doing so, we need to establish some standard definitions and results relating to weakly norming graphs and domination inequalities.
\subsection{Weakly norming graphs}

We begin by formally defining weakly norming graphs. Although weakly norming graphs are more commonly studied in the context of graph limits~\cite{GarbeHladkyLee22,ConlonLee17,KralMartinsPachWrochna19,LeeSchulke21}, we will instead follow~\cite[Definition~2.2]{Hatami10} and work in the equivalent setting of matrices in order to make the arguments in this paper more intuitive to combinatorialists.

Given an $n \times n$ matrix $W$ and $1\leq i,j\leq n$, let $W(i,j)$ denote the entry in the $i$th row and $j$th column of $W$. The definition of weakly norming graphs uses the following matrix generalisations of the homomorphism counting function. 

\begin{defn}
\label{def:weightedOriented}
Let $\vec{H}$ be an oriented graph. For $n\geq1$ and a matrix $W\in\mathbb{R}^{n\times n}$, define
\[h_{\vec{H}}(W) \coloneq \hom(\vec{H},W)\coloneq\sum_{\phi:V(\vec{H})\to[n]}\left(\prod_{(u,v)\in A(\vec{H})}W(\phi(u),\phi(v))\right).\]
Similarly, for an undirected graph $H$ and a symmetric matrix $W\in\mathbb{R}^{n\times n}$, define
\[h_H(W) \coloneq \hom(H,W)\coloneq\sum_{\phi:V(H)\to[n]}\left(\prod_{uv\in E(H)}W(\phi(u),\phi(v))\right).\]
\end{defn}

Observe that if $\vec{H}$ and $\vec{G}$ are oriented graphs and $W$ is the adjacency matrix of $\vec{G}$, then ${\hom(\vec{H},\vec{G})=\hom(\vec{H},W)}$; the analogous statement holds when $H$ is a graph and $W$ is the (symmetric) adjacency matrix of an undirected graph. Thus, the above definitions can be thought of as ``edge/arc-weighted'' versions of $\hom(\cdot,\cdot)$. Notice that when $W$ is symmetric, we have $h_{\vec{H}}(W) = h_H(W)$. We can now define weakly norming graphs.

\begin{defn}[Weakly norming]
    Let $H$ be a graph with at least one edge. For $n\geq1$ and a symmetric matrix $W\in \mathbb{R}^{n\times n}$, define
    \[\|W\|_{r(H)}\coloneq \hom(H,|W|)^{1/e(H)},\]
    where $|W|$ is the matrix obtained from $W$ by replacing each entry with its absolute value. A connected graph $H$ is \emph{weakly norming} if and only if $\|\cdot\|_{r(H)}$ is a norm on the set of real symmetric\footnote{Note that the definition of weakly norming graphs in~\cite[Definition~2.2]{Hatami10} does not require that the matrices are symmetric. However, it is shown in~\cite[Theorem~1.3]{LeeSidorenko22} that when $H$ is connected, it suffices to consider symmetric matrices.} $n\times n$ matrices for all $n\geq1$.
\end{defn}

\begin{rem}
\label{rem:triangleIneq}
It is clear that $\|cW\|_{r(H)}=|c|\cdot \|W\|_{r(H)}$ for any $W$ and scalar $c$, and so $\|\cdot\|_{r(H)}$ is homogeneous. Furthermore, $\|\cdot\|_{r(H)}$ is positive-definite if and only if $H$ is bipartite. Indeed, if $H$ is not bipartite and $W$ is the adjacency matrix of $K_2$, then $W \neq 0$ but $\hom(H,|W|)=0$, and so $\|\cdot\|_{r(H)}$ is not positive definite; on the other hand, if $W$ has a non-zero entry $W(i,j)$, then by mapping all vertices on one side of the bipartition to $i$ and the other side to $j$, we have that $\hom(H,|W|)>0$. Thus if $H$ is bipartite, showing that $\|\cdot\|_{r(H)}$ satisfies the triangle inequality is sufficient to prove that $H$ is weakly norming.
\end{rem}

\subsection{Domination inequalities}
We will also make use of an equivalent definition of relatively Sidorenko in terms of so-called `domination inequalities'. We say that an oriented graph $\vec{H}$ \emph{dominates} an oriented graph $\vec{F}$ if $\hom(\vec{F},\vec{G})\leq \hom(\vec{H},\vec{G})$ for every oriented graph $\vec{G}$. The proof will use a variant of the ``tensor power trick'', which is frequently used in the study of Sidorenko's conjecture~\cite{Sidorenko91,Hatami10,ConlonFoxSudakov10}.

\begin{figure}[htb]
\centering
\begin{tikzpicture}[scale=0.75, font=\small, x=1.6cm, y=1.3cm, >={Stealth[scale=0.7]}]

\node at (2.5,4.6) {$\vec{G}_2$};
\node at (-0.75,2.0) {$\vec{G}_1$};
\node at (2.5,0) {$\vec G_1\times \vec G_2$};

\node[mycircle,label=above:{$1$}] (h1) at (1,4) {};
\node[mycircle,label=above:{$2$}] (h2) at (2,4) {};
\node[mycircle,label=above:{$3$}] (h3) at (3,4) {};
\node[mycircle,label=above:{$4$}] (h4) at (4,4) {};

\draw[->] (h1) -- (h2);
\draw[->] (h2) -- (h3);
\draw[->] (h3) -- (h4);
\draw[->,bend left=25] (h3) to (h1);

\node[mycircle,label=left:{$a$}] (g1) at (0,3) {};
\node[mycircle,label=left:{$b$}] (g2) at (0,2) {};
\node[mycircle,label=left:{$c$}] (g3) at (0,1) {};

\draw[->] (g1) -- (g2);
\draw[->] (g2) -- (g3);

\begin{scope}[label distance=-2.5pt]
    \node[mycircle,label=above left:{$(a,1)$}]  (a1) at (1,3) {};
\node[mycircle,label=above:{$(a,2)$}]       (a2) at (2,3) {};
\node[mycircle,label=above:{$(a,3)$}]       (a3) at (3,3) {};
\node[mycircle,label=above right:{$(a,4)$}] (a4) at (4,3) {};

\node[mycircle,label=left:{$(b,1)$}]  (b1) at (1,2) {};
\node[mycircle]                       (b2) at (2,2) {};
\node[mycircle]                       (b3) at (3,2) {};
\node[mycircle,label=right:{$(b,4)$}] (b4) at (4,2) {};

\node[mycircle,label=below left:{$(c,1)$}]  (c1) at (1,1) {};
\node[mycircle,label=below:{$(c,2)$}]       (c2) at (2,1) {};
\node[mycircle,label=below:{$(c,3)$}]       (c3) at (3,1) {};
\node[mycircle,label=below right:{$(c,4)$}] (c4) at (4,1) {};
\end{scope}


\draw[->] (a1) -- (b2);
\draw[->] (a2) -- (b3);
\draw[->] (a3) -- (b1);
\draw[->] (a3) -- (b4);

\draw[->] (b1) -- (c2);
\draw[->] (b2) -- (c3);
\draw[->] (b3) -- (c1);
\draw[->] (b3) -- (c4);

\end{tikzpicture}
\caption{The tensor product of two oriented graphs.}
\label{fig:tensor-product-oriented}
\end{figure}
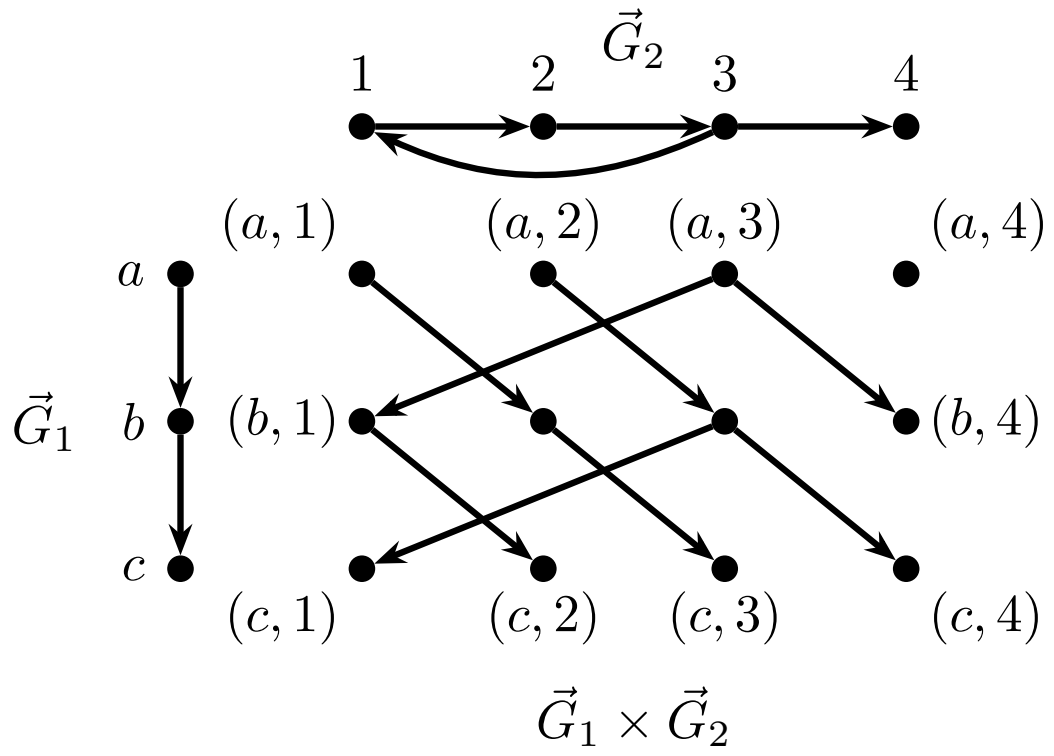

Recall that the \emph{tensor product} of oriented graphs $\vec{G}_1$ and $\vec{G}_2$ is the graph $\vec{G}_1\times\vec{G}_2$ with vertex set $V(\vec{G}_1)\times V(\vec{G}_2)$ which contains an arc from $(u_1,u_2)\in V(\vec{G}_1)\times V(\vec{G}_2)$ to $(v_1,v_2)\in V(\vec{G}_1)\times V(\vec{G}_2)$ if and only if $(u_1,v_1)\in A(\vec{G}_1)$ and $(u_2,v_2)\in A(\vec{G}_2)$. An example is given in \cref{fig:tensor-product-oriented}.
It is well-known and easy to verify that a homomorphism from $\vec{H}$ to a tensor product $\vec{G}_1\times \vec{G}_2$ is equivalent to a pair of homomorphisms from $\vec{H}$ to $\vec{G}_1$ and to $\vec{G}_2$.
\begin{restatable}[Folklore; see, e.g.,~{\cite[Corollary~2.3]{HellNesetril04}}]{lem}{tensorOrient}
\label{lem:tensorOrient}
For any oriented graphs $\vec{H}, \vec{G_1}$ and $\vec{G}_2$,
\[\hom(\vec{H},\vec{G}_1\times \vec{G}_2)=\hom(\vec{H},\vec{G}_1)\cdot\hom(\vec{H},\vec{G}_2).\]
\end{restatable}

The following result is an equivalent definition for relatively Sidorenko graphs, in terms of a domination property.

\begin{lem}
\label{lem:relSidDom}
Let $\vec{H}$ be an oriented graph. Then $\vec{H}$ is relatively Sidorenko if and only if it dominates every other orientation of its underlying graph.
\end{lem}

\begin{proof}
Let $\vec{H}_1,\vec{H}_2,\dots,\vec{H}_{2^{e(H)}}$ denote the distinct orientations of $H$ (not identifying isomorphic copies). Suppose first that $\vec{H}$ dominates $\vec{H}_i$ for all $1\leq i\leq 2^{e(H)}$. Then, for any oriented graph $\vec{G}$, we have
\[\hom(H,G)=\sum_{i=1}^{2^{e(H)}}\hom(\vec{H}_i,\vec{G}) \leq \sum_{i=1}^{2^{e(H)}}\hom(\vec{H},\vec{G}) = 2^{e(H)}\hom(\vec{H},\vec{G})\]
and so $\vec{H}$ is relatively Sidorenko. Conversely, suppose that there exists an orientation of $H$, say $\vec{H}_1$, and an oriented graph $\vec{G}$ such that $\hom(\vec{H}_1,\vec{G})>\hom(\vec{H},\vec{G})$. Let $r$ be a positive integer chosen sufficiently large so that
\begin{equation} \label{eq:amplify}
    \hom(\vec{H}_1,\vec{G})^r>2^{e(H)}\hom(\vec{H},\vec{G})^r.
\end{equation}
Let 
\[\vec{G}^{\times r}=\underbrace{\vec G\times\cdots\times\vec G}_{r}\] 
be the $r$-fold tensor product of $\vec{G}$ with itself and let $G'$ be the underlying graph of $\vec{G}^{\times r}$. Observe that  $\hom(H,G')\geq \hom(\vec{H}_1,\vec{G}^{\times r})$ since every homomorphism from $\vec{H}_1$ to $\vec{G}^{\times r}$ yields a homomorphism from $H$ to $G'$. Therefore, by \cref{lem:tensorOrient} and \eqref{eq:amplify}, we have
\[\hom(H,G')\geq \hom(\vec{H}_1,\vec{G}^{\times r}) = \hom(\vec{H}_1,\vec{G})^r>2^{e(H)}\hom(\vec{H},\vec{G})^r=2^{e(H)}\hom(\vec{H},\vec{G}^{\times r})\]
and so $\vec{H}$ is not relatively Sidorenko.  
\end{proof}

\subsection{Proofs of \cref{th:relSidCharacterisation} and \cref{th:relSidExamples}}

We now have the background required to prove \cref{th:relSidCharacterisation}. We begin by establishing three necessary conditions for an oriented graph to be relatively Sidorenko.

\begin{lem}
\label{lem:relSidBip}
Every relatively Sidorenko oriented graph is bipartite.
\end{lem}
\begin{proof}
Let $\vec{H}$ be a relatively Sidorenko graph, and suppose for a contradiction that $H$ is not bipartite. We will construct an oriented graph $\vec{G}$ such that $\hom(\vec{H},\vec{G}) = 0$ but ${\hom(H,G) > 0}$, contradicting the relative Sidorenko definition.

Let $C$ be the \emph{core} of $H$. That is, $C$ is the vertex-minimal subgraph of $H$ such that there exists a homomorphism $\rho:H\to C$. Clearly, the restriction of $\rho$ to $C$ is a graph homomorphism from $C$ to itself. Furthermore, $\rho^2 \colon H \to H$ is also a graph homomorphism, and its image is ${\rho^2(V(H)) = \rho(V(C)) = V(C)}$; otherwise, $\rho^2$ would contradict the minimality of $C$. Thus $\rho|_C \colon C \to C$ is a bijection, and therefore an isomorphism from $C$ to itself. Furthermore, any graph homomorphism $\sigma : C \to C$ must be an isomorphism from $C$ to itself, otherwise $\sigma \circ \rho$ would also contradict the choice of $\rho$.

Since $H$ is not bipartite, it further follows that $C$ is not bipartite. Indeed, if $C$ were bipartite then there would be a graph homomorphism $\sigma \colon C \to K_2$, and composing $\sigma$ with $\rho \colon H \to C$ would construct a graph homomorphism $H \to K_2$. Since $H$ is not bipartite, such a homomorphism does not exist. It follows that $C$ contains an odd cycle. 

Let $\vec{C_0}$ be the orientation of $C$ induced by $\vec{H}$. We choose another orientation $\vec{C_1}$ of the same underlying graph $C$ such that $\vec{C_1}$ is not obtained from $\vec{C_0}$ by any isomorphism $C \to C$. If $\vec{C_0}$ is acyclic, construct $\vec{C_1}$ by orienting the edges of some odd cycle in $C$ cyclically, and choosing an arbitrary orientation for all remaining edges of $C$. Since $\vec{C_1}$ contains a directed cycle, it is necessarily not the image of the acyclic $\vec{C_0}$ under any isomorphism $C \to C$. Otherwise, if $\vec{C_0}$ does contain a directed cycle, choose $\vec{C_1}$ to be any acyclic orientation of $C$. It follows similarly that $\vec{C_1}$ is an orientation of $C$ that is not obtained by an isomorphism $C \to C$ applied to $\vec{C_0}$.

Now define $\vec{G} \coloneq \vec{C_1}$ and $G \coloneq C$. Since $\rho \colon H \to C=G$ is a graph homomorphism, we have $\hom(H,G) > 0$. Suppose, for a contradiction, that there is also a homomorphism $\varphi \colon \vec{H} \to \vec{C_1}$. Then since $\varphi$ is a graph homomorphism $H \to C$, it follows that the restriction $\varphi|_C \colon C \to C$ is an isomorphism from $C$ to itself, and so $\vec{C_1}$ is precisely the image of $\vec{C_0}$ under $\varphi$. But this contradicts the construction of $\vec{C_1}$, and so $\varphi$ cannot exist. It follows that $\hom(\vec{H}, \vec{C_1}) = 0$, contradicting the assumption that $\vec{H}$ is relatively Sidorenko.
\end{proof}

\begin{prop}
\label{prop:relSidOneWay}
Let $H$ be a graph. Then every relatively Sidorenko orientation of $H$ is one-way and has the same number of sources as sinks. 
\end{prop}

\begin{proof}
We may assume $H$ is bipartite; otherwise, by \cref{lem:relSidBip}, it has no relatively Sidorenko orientations and the statement holds vacuously. Suppose that $\vec{H}$ is an orientation of $H$ and let $\vec{K}_2$ be the unique oriented graph with two vertices and one arc. Since $H$ is bipartite, $\hom(H,K_2)>0$. If the orientation of $\vec{H}$ is not one-way, then $\hom(\vec{H},\vec{K}_2)=0$ and so $\vec{H}$ is not relatively Sidorenko.

Now, suppose that the orientation of $\vec{H}$ is one-way and that $\vec{H}$ has $a$ sources and $b < a$ sinks. Let $c$ be the number of isolated vertices of $\vec{H}$ (which are both sources and sinks). Let $\vec{G}$ be an oriented star with $n\geq2$ leaves, with all arcs directed away from the centre. Then $\hom(\vec{H},\vec{G})=n^{b-c}(n+1)^c$ since every source that is not isolated must map to the central vertex of $\vec{G}$. On the other hand, since we may choose which vertex part of $H$ is mapped to the central vertex of $G$, we have $\hom(H,G)\geq \left(n^{a-c} + n^{b-c}\right)(n+1)^c \geq n^{a-c}(n+1)^c$. Thus, $\hom(H,G)\geq n^{a-b}\hom(\vec{H},\vec{G})$. Since $a>b$, it follows that $\hom(H,G)$ is greater than $2^{e(H)}\hom(\vec{H},\vec{G})$ for all sufficiently large $n$. Therefore, $\vec{H}$ is not relatively Sidorenko. The case where $a < b$ is similar, completing the proof.
\end{proof}

\begin{lem}
\label{lem:relSidReversible}
 If $\vec{H}$ is a connected relatively Sidorenko oriented graph with at least one arc, then $H$ is reversible.
\end{lem}

\begin{proof}
Recall that $\vec{H}^{\text{rev}}$ is the oriented graph obtained from $\vec{H}$ by reversing the orientations of all of its arcs. Suppose that $\vec{H}$ is a relatively Sidorenko oriented graph. By \cref{lem:relSidBip} and \cref{prop:relSidOneWay}, we know that $H$ is bipartite and $\vec{H}$ is a one-way orientation. Consider an arbitrary oriented host graph $\vec{G}$. By \cref{lem:relSidDom}, since $\vec{H}$ is relatively Sidorenko, it dominates all other orientations of $H$, and so $\hom(\vec{H},\vec{G}) \geq \hom(\vec{H}^{\text{rev}}, \vec{G})$. Furthermore, by considering the effect of reversing the arc directions in $\vec{H}$ and $\vec{G}$ simultaneously, we have that
\begin{align*}
    \hom(\vec{H}, \vec{G}) = \hom(\vec{H}^{\text{rev}}, \vec{G}^{\text{rev}}) \leq \hom(\vec{H}, \vec{G}^{\text{rev}}) = \hom(\vec{H}^{\text{rev}}, \vec{G}),
\end{align*}
and so it follows that $\hom(\vec{H}, \vec{G}) = \hom(\vec{H}^{\text{rev}},\vec{G})$ for all oriented host graphs $\vec{G}$. Now, for any oriented graph $\vec{F}$, let $\text{sur}(\vec{F},\vec{H})$ denote the number of homomorphisms from $\vec{F}$ to $\vec{H}$ that map to every vertex of $\vec{H}$. By standard inclusion-exclusion, we have
\begin{equation*}
    \text{sur}(\vec{F},\vec{H}) = \sum_{S \subseteq V(\vec{H})}(-1)^{|V(\vec{H})|-|S|}\hom(\vec{F},\vec{H}[S]).
\end{equation*}
By the previous argument, we have that $\hom(\vec{H},\vec{H}[S]) = \hom(\vec{H}^{\text{rev}},\vec{H}[S])$ for all $S \subseteq V(\vec{H})$, and so it follows that $\text{sur}(\vec{H}^{\text{rev}}, \vec{H}) = \text{sur}(\vec{H}, \vec{H}) \geq 1$. Since any homomorphism from $\vec{H}^{\text{rev}}$ to $\vec{H}$ that maps to every vertex of $\vec{H}$ must be an automorphism, this completes the proof that $H$ is reversible. 
\end{proof}

We now prove our main characterisation of relatively Sidorenko oriented graphs in terms of weakly norming graphs, and use this characterisation to derive examples of the former.

\relSidCharacterisation*

\begin{proof}
By \cref{lem:relSidBip,lem:relSidReversible}, \cref{prop:relSidOneWay} and the hypotheses of the theorem, we can assume throughout the proof that $H$ is connected, bipartite, reversible and contains at least one edge, and that $\vec{H}$ is a one-way orientation of $H$. Let $(X,Y)$ be the bipartition of $H$ so that all arcs of $\vec{H}$ are directed from $X$ to $Y$. Under the aforementioned conditions, we show that $\vec{H}$ is relatively Sidorenko if and only if $H$ is weakly norming.

For simplicity of notation, we let $m = e(H)$, let $v = v(H)$, and denote both $h_H$ and $h_{\vec{H}}$ simply by $h$ throughout this proof. This will cause no ambiguity, as we only work with symmetric matrices and $h_{\vec{H}}(W) = h_H(W)$ whenever $W$ is symmetric. We begin by showing that if $H$ is a reversible weakly norming graph, then the one-way orientation $\vec{H}$ of $H$ is relatively Sidorenko.
\begin{claim}
    Suppose that $H$ is a connected, reversible weakly norming graph with at least one edge. Then the one-way orientation $\vec{H}$ is relatively Sidorenko.
\end{claim}
\begin{poc}
    We will construct a host graph that allows us to deduce the relatively Sidorenko definition from the norm triangle inequality on $h(|\cdot|)^{1/m}$. The definition of relatively Sidorenko holds trivially for homomorphisms from $\vec{H}$ into the empty host graph, so consider a non-empty host $\vec{G}$. From $\vec{G}$, construct the auxiliary bipartite graph $B$, whose vertex parts $V(G)_L$ and $V(G)_R$ are both copies of $V(G)$. For the copy $x_L$ of $x$ in $V(G)_L$ and $y_R$ of $y$ in $V(G)_R$, put the edge $x_Ly_R$ in $B$ if and only if $(x,y)$ is an arc of $\vec{G}$. Let $B'$ be the graph obtained from $B$ by swapping its two vertex parts such that $x_L \leftrightarrow x_R$ for each $x \in V(G)$, and note that $B \cong B'$. Furthermore, since an edge common to both $B$ and $B'$ would require both arc directions between a pair of vertices in $\vec{G}$, the two graphs are edge-disjoint. The union graph $B \cup B'$ on their common vertex set contains precisely the edges $x_Ly_R$ and $y_Lx_R$ for each $xy$ in $E(G)$. 
    
    Since $H$ is connected, any homomorphism of $H$ into $B$ maps all vertices of $X$ to one vertex part of $B$, and all vertices of $Y$ to the other. A homomorphism mapping $X$ to $V(G)_L$ corresponds to a homomorphism from $\vec{H}$ to $\vec{G}$, and the opposite placement corresponds to a homomorphism from $\vec{H}^{\text{rev}}$ to $\vec{G}$. Since $H$ is reversible, the number of homomorphisms of each type is equal, and so
    
    \begin{equation}\label{eq:homHB}
        \hom(H,B) = \hom(\vec{H},\vec{G}) + \hom(\vec{H}^{\text{rev}},\vec{G}) = 2\hom(\vec{H},\vec{G}).
    \end{equation}

    By similar reasoning, every homomorphism from $H$ to $B \cup B'$ can be realised by either mapping $X$ to $V(G)_L$ and $Y$ to $V(G)_R$, or vice versa. Since both cross edges $x_Ly_R$ and $y_Lx_R$ are present in $B \cup B'$, the edges of $G$ are present in either placement. Hence, $\hom(H, B \cup B') = 2\hom(H,G)$.  

    Finally, note that the adjacency matrices of $B$ and $B'$ are symmetric matrices on the same vertex set, and the sum of these matrices is the adjacency matrix of $B \cup B'$. Since $H$ is weakly norming, the triangle inequality gives
    \begin{align*}
        \hom(H, B \cup B')^{1/m} \leq \hom(H,B)^{1/m} + \hom(H,B')^{1/m} &= 2\hom(H,B)^{1/m},
     \end{align*}
     and hence
    \begin{align*}
        2\hom(H,G) = \hom(H, B \cup B') &\leq 2^{m}\hom(H,B).
    \end{align*}
    Substituting this back into \eqref{eq:homHB} completes the proof.
\end{poc}

Suppose now that $\vec{H}$ is relatively Sidorenko. By \cref{rem:triangleIneq}, to show that $H$ is weakly norming, it suffices to show that $h(|\cdot|)^{1/m}$ satisfies the norm triangle inequality. We first prove that $h$ satisfies the midpoint inequality on real, non-negative, symmetric matrices. This will in turn show that $h$ is convex on this space, from which we will derive the triangle inequality.
\begin{claim}
\label{cl:midpointConvex}
For any symmetric non-negative $n\times n$ real matrices $A$ and $B$, it holds that
\[h\left(\frac{A+B}{2}\right)\leq \frac{h(A)+h(B)}{2}.\]
\end{claim}

\begin{poc}
Suppose for now that $A$ and $B$ satisfy $A(i,j) + B(i,j) \leq 1$ for all pairs $i,j \in [n]$. We construct an oriented graph in which these two weights encode the probabilities of two possible directions of an edge.

Construct a random bipartite oriented graph $\vec{G}_N$ as follows. Take disjoint vertex classes $L_1,\dots,L_n$ and $R_1,\dots,R_n$, each of size $N$. Between each pair $x \in L_i$, $y \in R_j$, independently take the arc $(x,y)$ with probability $A(i,j)$, the arc $(y,x)$ with probability $B(i,j)$, and no edge with probability $1 - A(i,j) - B(i,j)$. We will count $\hom(\vec{H},\vec{G}_N)$, and choose $N$ to be sufficiently large so that the contribution of non-injective homomorphisms to this sum is negligible. 

Since $\vec{H}$ is relatively Sidorenko, we have that $H$ is a connected bipartite graph and $\vec{H}$ is the one-way orientation of $H$. Let $X$ and $Y$ denote the source and sink vertex classes of $\vec{H}$. Since $H$ is connected, either $X$ maps into $L_1 \cup \dots \cup L_n$ and $Y$ maps into $R_1\cup\dots\cup R_n$, or vice versa. Suppose that $X$ maps to $L_1 \cup \dots \cup L_n$, and fix any map $\phi : V(H) \to [n]$. A \emph{lift} of $\phi$ is a map from $V(\vec{H})$ to $V(\vec{G}_N)$ such that each $x \in X$ is mapped to a vertex in $L_{\phi(x)}$, and each vertex $y \in Y$ is mapped to a vertex in $R_{\phi(y)}$. For any choice of $\phi$, there are $N^{v}$ such lifts, which may or may not be injective. For any pair of vertices in $\vec{H}$, at most $N^{v - 1}$ of the lifts identify that pair, so the total number of non-injective homomorphisms from $\vec{H}$ to $\vec{G}_N$ is $\mathcal{O}_{n,v}(N^{v-1})$ (we consider $n$ and $v$ to be constants). 

Whenever a lift is injective, the arcs of $\vec{H}$ exist independently in $\vec{G}_N$, and so any given injective lift gives an injective homomorphism between $\vec{H}$ and $\vec{G}_N$ with probability
$$ \prod_{xy \in E(H)} A(\phi(x), \phi(y)).$$ Similar calculations apply in the case where $X$ is mapped into $R_1\cup\dots\cup R_n$, so the sum of the probabilities that a lift of $\phi$ is an injective homomorphism over both possible vertex placements is
$$p_\phi \coloneq \prod_{xy \in E(H)} A(\phi(x), \phi(y)) + \prod_{xy \in E(H)} B(\phi(x), \phi(y)).$$
Since $p_\phi \leq 1$, it follows that the expected number of injective homomorphisms arising from lifts of a given map $\phi$, summed over both vertex placements is $N^vp_\phi + \mathcal{O}_{n,v}(N^{v-1})$. Hence, counting the total expected number of injective and non-injective homomorphisms between $\vec{H}$ and $\vec{G}_N$, we have
\begin{align*}
    \mathbb{E}[\hom(\vec{H},\vec{G}_N)] &= \sum_{\phi: V(H) \to [n]}N^vp_\phi + \mathcal{O}_{n,v}(N^{v-1})\\
    &= N^{v}(h(A) + h(B)) + \mathcal{O}_{n,v}(N^{v-1}).
\end{align*}
It follows that as $N \to \infty$, we have 
\begin{equation}
    \frac{1}{N^{v}}\mathbb{E}[\hom(\vec{H},\vec{G}_N)] \to h(A) + h(B).   
\end{equation}
Applying the same strategy to count $\hom(H, G_N)$, we find that the probability of realising a given edge $xy \in E(H)$ where $\phi(x) = i$ and $\phi(y) = j$ is $A(i,j) + B(i,j)$ when $x$ maps to $L_i$ and $y$ to $R_j$, and by symmetry of $A$ and $B$, the same is true when $x$ is mapped into $R_i$ and $y$ into $L_j$. Hence, in the undirected case, the sum over both possible vertex placements of the probabilities that a given lift of some $\phi : V(H) \to [n]$ gives an injective homomorphism from $H$ to $G_N$ is 
$$q_\phi \coloneq 2\prod_{xy \in E(H)}(A+B)(\phi(x), \phi(y)).$$

Following analogous calculations to the directed case, we therefore find that as $N \to \infty$,
\begin{equation}
    \frac{1}{N^{v}}\mathbb{E}[\hom(H,G_N)] \to 2h(A+B).  
\end{equation}
Since $\vec{H}$ is relatively Sidorenko, it follows by definition that 
\begin{align*}
    \hom(\vec{H},\vec{G}_N) &\geq 2^{-m}\hom(H,G_N)\\
    \implies \mathbb{E}[\hom(\vec{H},\vec{G}_N)] &\geq 2^{-m}\mathbb{E}[\hom(H,G_N)]\\
    \implies h(A) + h(B) &\geq \frac{1}{2^{m-1}}h(A+B)\\
    \implies \frac{h(A) + h(B)}{2} &\geq h\left(\frac{A+B}{2}\right),
\end{align*}
and so the midpoint inequality holds for any symmetric non-negative $n \times n$ matrices $A$ and $B$ satisfying $A(i,j) + B(i,j) \leq 1$ for all $i$,$j$. Finally, noting that $h(cW) = c^{m}h(W)$ for any non-negative symmetric $n \times n$ matrix $W$, the result follows for arbitrary symmetric non-negative $n \times n$ matrices by applying an appropriate scaling factor. This completes the proof. 
\end{poc}

We now show that this implies that $h$ is convex on the space of symmetric non-negative matrices.
\begin{claim}\label{clm:convex}
    For any symmetric non-negative $n \times n$ real matrices $A$ and $B$, and any $\lambda \in [0,1]$, we have that
    $$h((1-\lambda)A + \lambda B) \leq (1- \lambda)h(A) + \lambda h(B).$$
\end{claim}
\begin{poc}
    \cref{cl:midpointConvex} shows that this holds when $\lambda = \frac{1}{2}$. By repeatedly applying \cref{cl:midpointConvex} to obtain appropriate midpoints, we can extend this to show that the claim holds for any dyadic $\lambda \in [0,1]$. Finally, observe that for any fixed $A$ and $B$, the map $\lambda \to h((1-\lambda)A + \lambda B)$ is a polynomial, and hence continuous. Thus, since the dyadic rationals are dense in $[0,1]$, the statement holds for all $\lambda \in [0,1]$, as desired.
\end{poc}
We complete this direction of the proof by showing that $h(|\cdot|)^{1/m}$ satisfies the triangle inequality, and is therefore a norm on the space of real symmetric $n \times n$ matrices. 
\begin{claim}
    Let $A$ and $B$ be  real symmetric $n \times n$ matrices. Then 
    $$h(|A+B|)^{1/m} \leq h(|A|)^{1/m} + h(|B|)^{1/m}.$$
\end{claim}
\begin{poc}
    Begin by assuming that $A$ and $B$ are both non-negative matrices. If either is the zero matrix, we are clearly done, so assume that both contain at least one non-zero entry. Note that by definition of $h$, we have that $$h\left(\frac{X}{h(X)^{1/m}}\right) = \frac{h(X)}{(h(X)^{1/m})^{m}} = 1,$$
    for any non-negative, non-zero symmetric matrix $X$. Letting ${\lambda = \frac{h(B)^{1/m}}{h(A)^{1/m} + h(B)^{1/m}}}$ and applying \cref{clm:convex}, we have that
    \begin{align*}
        \frac{h(A+B)}{(h(A)^{1/m} + h(B)^{1/m})^m} &= h\left(\frac{A+B}{h(A)^{1/m} + h(B)^{1/m}}\right)\\
        &= h\left(\frac{h(A)^{1/m}}{h(A)^{1/m} + h(B)^{1/m}}\cdot \frac{A}{h(A)^{1/m}} + \frac{h(B)^{1/m}}{h(A)^{1/m} + h(B)^{1/m}}\cdot \frac{B}{h(B)^{1/m}}\right) \\
        &\leq \frac{h(A)^{1/m}}{h(A)^{1/m} + h(B)^{1/m}}\cdot h\left(\frac{A}{h(A)^{1/m}}\right) + \frac{h(B)^{1/m}}{h(A)^{1/m} + h(B)^{1/m}}\cdot h\left(\frac{B}{h(B)^{1/m}}\right) \\
        &= \frac{h(A)^{1/m}}{h(A)^{1/m} + h(B)^{1/m}} + \frac{h(B)^{1/m}}{h(A)^{1/m} + h(B)^{1/m}}\\
        &= 1,
    \end{align*}
    and so $h(A+B)^{1/m} \leq h(A)^{1/m} + h(B)^{1/m}$. For arbitrary real symmetric matrices $A$ and $B$, observe that we have $|A + B| \leq |A| + |B|$ entry-wise, and that $h$ is entry-wise increasing on non-negative matrices. Hence, the claim follows for all real symmetric matrices by applying the previous calculation to $|A|$ and $|B|$.
\end{poc}
It now follows from \cref{rem:triangleIneq} that $H$ is weakly norming, completing the proof.
\end{proof}

While a characterisation of weakly norming graphs entirely in terms of graph structure remains unknown, strong necessary conditions are known. For example, Sidorenko~\cite{Sidorenko20} proved that every weakly norming graph is edge-transitive. While weakly norming graphs in general need not be vertex-transitive (for example, complete bipartite graphs $K_{n,m}$ with $n \neq m$ are weakly norming), vertex transitivity follows easily from the combination of edge-transitivity and reversibility, establishing the following corollary.
\begin{cor}
    The underlying graph of every connected relatively Sidorenko graph with at least one arc is both vertex- and edge-transitive.
\end{cor}

We conclude by listing some explicit families of relatively Sidorenko graphs.
\relSidExamples*

\begin{proof}
 It is easily observed that the graphs in each of these classes are reversible. Note that family (c) is disconnected for $k =1,2$. However, both these graphs are easily seen to be relatively Sidorenko by checking the definition. In all other cases, the families consist of connected, reversible bipartite graphs containing at least one arc. Now, even cycles are weakly norming due to a connection to Schatten--von Neumann norms observed by Hatami~\cite[Section~2.3]{Hatami10}. Hatami~\cite[Theorem~2.9]{Hatami10} also proved that complete bipartite graphs and hypercubes are weakly norming, and Lov\'asz~\cite[Proposition~14.2(c)]{Lovasz12} proved that balanced complete bipartite graphs with a perfect matching removed are weakly norming. Thus, the result follows by \cref{th:relSidCharacterisation}.
\end{proof}

\section{Relatively Anti-Sidorenko Oriented Graphs}
\label{sec:antiSid}
\begin{figure}[htbp]
    \centering
    \begin{tikzpicture}[scale=0.8,
  thick,
  >={Stealth[scale=1.2]},
  dot/.style={circle, fill=black, inner sep=2pt}
]
\begin{scope}[scale=0.7]

\node[dot] (r) at (0,0) {};

\node[dot] (L) at (-2.5,-1.8) {};
\node[dot] (R) at ( 2.5,-1.8) {};

\draw[->] (r) -- (L);
\draw[->] (R) -- (r);

\foreach \j/\x in {
  1/-4.375,
  2/-3.125,
  3/-1.875,
  4/-0.625
}{
  \node[dot] (L2\j) at (\x,-3.6) {};
}

\foreach \j/\x in {
  1/0.625,
  2/1.875,
  3/3.125,
  4/4.375
}{
  \node[dot] (R2\j) at (\x,-3.6) {};
}

\draw[->] (L) -- (L21);
\draw[->] (L) -- (L22);
\draw[->] (L23) -- (L);
\draw[->] (L24) -- (L);

\draw[->] (R) -- (R21);
\draw[->] (R) -- (R22);
\draw[->] (R23) -- (R);
\draw[->] (R24) -- (R);

\def\dx{0.29}

\foreach \j in {1,2,3,4}{
  \node[dot] (L3a\j)
    at ($(L2\j)+(-\dx,-1.8)$) {};
  \node[dot] (L3b\j)
    at ($(L2\j)+(\dx,-1.8)$) {};

  \draw[->] (L2\j) -- (L3a\j);
  \draw[->] (L3b\j) -- (L2\j);
}

\foreach \j in {1,2,3,4}{
  \node[dot] (R3a\j)
    at ($(R2\j)+(-\dx,-1.8)$) {};
  \node[dot] (R3b\j)
    at ($(R2\j)+(\dx,-1.8)$) {};

  \draw[->] (R2\j) -- (R3a\j);
  \draw[->] (R3b\j) -- (R2\j);
}

\end{scope}

\begin{scope}[scale=0.7,xshift=12cm]

\node[dot] (u) at (-3,0) {};
\node[dot] (v) at ( 3,0) {};
\draw[->] (u) -- (v);

\node[dot] (uout) at (-4.2,-2.7) {};
\node[dot] (uin)  at (-1.8,-2.7) {};

\node[dot] (vin)  at ( 1.8,-2.7) {};
\node[dot] (vout) at ( 4.2,-2.7) {};

\draw[->] (u) -- (uout);
\draw[->] (uin) -- (u);

\draw[->] (vin) -- (v);
\draw[->] (v) -- (vout);

\node[dot] (uoutL) at (-4.8,-5.4) {};
\node[dot] (uoutR) at (-3.6,-5.4) {};
\draw[->] (uout) -- (uoutL);
\draw[->] (uoutR) -- (uout);

\node[dot] (uinL) at (-2.4,-5.4) {};
\node[dot] (uinR) at (-1.2,-5.4) {};
\draw[->] (uin) -- (uinL);
\draw[->] (uinR) -- (uin);

\node[dot] (vinL) at ( 1.2,-5.4) {};
\node[dot] (vinR) at ( 2.4,-5.4) {};
\draw[->] (vin) -- (vinL);
\draw[->] (vinR) -- (vin);

\node[dot] (voutL) at ( 3.6,-5.4) {};
\node[dot] (voutR) at ( 4.8,-5.4) {};
\draw[->] (vout) -- (voutL);
\draw[->] (voutR) -- (vout);

\end{scope}
\end{tikzpicture}
\caption{Examples of two perfectly balanced oriented trees. On the left, 
$|S_0|=1$, $k=3$, $d_1=1$, $d_2=2$ and $d_3=1$; on the right, 
$|S_0|=2$, $k=2$, $d_1=1$ and $d_2=1$.}
    \label{fig:balancedTrees}
\end{figure}
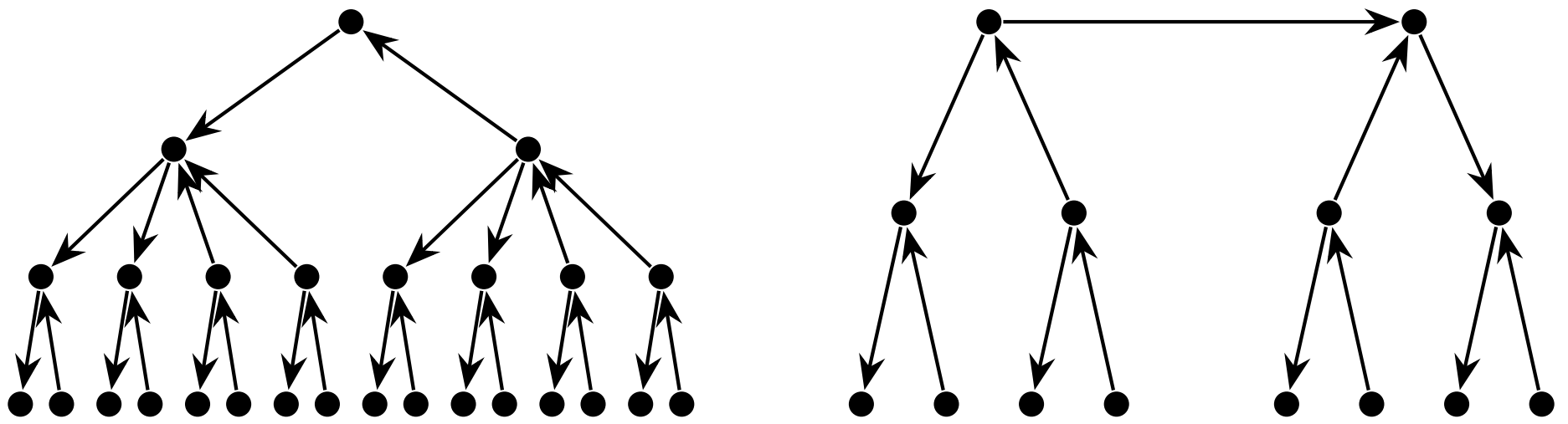

In this section we construct some families of relatively anti-Sidorenko oriented graphs. The first family is the family of \emph{perfectly balanced} trees, which we now define. Let $k\geq0$ and let $d_1,\dots,d_k$ be positive integers, and let $S_0$ be either a single vertex or two vertices joined by an arc. For each $1 \leq i \leq k$, attach $2d_i$ new leaf neighbours to each vertex $v\in S_{i-1}$, with $d_i$ of the arcs directed toward $v$ and $d_i$ arcs directed away from $v$, and let $S_i$ be the set of all such leaves. We say that such an oriented tree is perfectly balanced. Some examples of perfectly balanced oriented trees are given in \cref{fig:balancedTrees}.

\perfectlyBalanced*

In fact, we will prove the following stronger result.

\begin{thm}
\label{th:antiSidTreesStronger}
Let $\vec{T}$ be a perfectly balanced oriented tree and let $\vec{T}_1,\dots,\vec{T}_{2^{e(T)}}$ enumerate all orientations of $T$, 
not identifying isomorphic copies. Then, for any oriented graph $\vec{G}$,
\[\hom(\vec{T},\vec{G})\leq \left(\prod_{i=1}^{2^{e(T)}}\hom(\vec{T}_i,\vec{G})\right)^{1/2^{e(T)}}.\]
\end{thm}

The first part of this section is focused on the proof of \cref{th:antiSidTreesStronger}, from which \cref{th:antiSidTrees} will follow. We take an information-theoretic approach which originated in a paper of Kopparty and Rossman~\cite{KoppartyRossman11} and which is also used in, e.g.,~\cite{BlekhermanRaymond23,BehagueCrudeleNoelSimbaqueba25,ChenClemenNoel25+,BehagueMorrisonNoel24,ChenClemenNoel26+}. Some standard definitions and tools that we will need are presented in \cref{ssec:entropy}, and our argument is given in \cref{ssec:entropy_proof}. Finally, in \cref{ssec:non_bpt_antiSid}, we give a recursive construction of relatively anti-Sidorenko oriented graphs with arbitrary chromatic number, proving \cref{th:antiSidChrom}. This is in sharp contrast to the behaviour for relatively Sidorenko graphs, which by \cref{lem:relSidBip} have chromatic number at most $2$.

\subsection{Entropy preliminaries}
\label{ssec:entropy}

Throughout this section, all random variables have finite range. Let $X$ be a discrete random variable. The \emph{range} $\rng(X)$ of $X$ is the set of all $x$ such that $\mathbb{P}(X=x) > 0$. The \emph{entropy} of $X$ is 
\[\mathbb{H}(X) := \sum_{x\in\rng(X)}\mathbb{P}(X=x)\log_2\left(\frac{1}{\mathbb{P}(X=x)}\right).\]
The convexity of $x\log_2(x)$ implies that the entropy of $X$ is maximised when $X$ is uniformly distributed over $\rng(X)$. That is, we have the following standard result.
\begin{lem}[Maximum Entropy]
\label{lem:uniform}
Every random variable $X$ with finite range satisfies
\[\mathbb{H}(X)\leq \log_2(|\rng(X)|),\]
where equality holds if and only if $\mathbb{P}(X=x)=\frac{1}{|\rng(X)|}$ for all $x\in\rng(X)$. 
\end{lem}
We will also use the standard notion of conditional entropy. Let $X$ and $Y$ be discrete random variables, and let $y \in \rng(Y)$. We denote by $\rng(X \mid Y=y)$ the set of all $x$ such that ${\mathbb{P}(X=x \mid Y=y) > 0}$. The entropy of $X$ conditioned on the event $Y=y$ is
\[\mathbb{H}(X\mid Y=y):=\sum_{x\in \rng(X\mid Y=y)}\mathbb{P}(X=x\mid Y=y)\log_2\left(\frac{1}{\mathbb{P}(X=x\mid Y=y)}\right).\]
The \emph{conditional entropy} of $X$ \emph{given} $Y$ is then defined as
\[\mathbb{H}(X\mid Y):=\sum_{y\in \rng(Y)}\mathbb{P}(Y=y)\mathbb{H}(X\mid Y=y).\]
The following standard fact allows us to express the entropy of a tuple of random variables in terms of conditional entropy.
\begin{lem}[Chain Rule for Entropy]
\label{lem:ChainRule}
For any discrete random variables $X_1,\dots,X_k$,
\[\mathbb{H}(X_1,\dots,X_k)=\mathbb{H}(X_1)+\sum_{i=2}^k\mathbb{H}(X_i\mid X_1,\dots,X_{i-1}).\]
\end{lem}

When applying the Chain Rule, it will often be useful to make use of the notion of conditional independence to simplify the conditional entropy expressions. For discrete random variables $X,Y$ and $Z$, we say that $X$ and $Z$ are \emph{conditionally independent given $Y$} if, for all $y\in\rng(Y)$, $x\in \rng(X\mid Y=y)$ and $z\in\rng(Z\mid Y=y)$, it holds that
\[\mathbb{P}(X=x\mid Y=y,Z=z)=\mathbb{P}(X=x\mid Y=y),\]
and
\[\mathbb{P}(Z=z\mid Y=y,X=x)=\mathbb{P}(Z=z\mid Y=y).\]

\begin{lem}[Deconditioning]
\label{lem:decon}
For any three discrete random variables $X,Y$ and $Z$,
\[\mathbb{H}(X\mid Y,Z)\leq \mathbb{H}(X\mid Y)\]
where equality holds if and only if $X$ and $Z$ are conditionally independent given $Y$. 
\end{lem}

\subsection{Proof of \texorpdfstring{\cref{th:antiSidTreesStronger}}{Theorem 3.1}}
\label{ssec:entropy_proof}

We now derive an expression for the entropy of a uniformly random homomorphism from an oriented tree to an oriented graph. 
\begin{lem}
\label{lem:uniformHom}
Let $\vec{T}$ be an oriented tree and let $\vec{G}$ be an oriented graph such that ${\hom(\vec{T},\vec{G})>0}$. If $\varphi:V(\vec{T})\to V(\vec{G})$ is a uniformly random homomorphism from $\vec{T}$ to $\vec{G}$, then
\[\mathbb{H}(\varphi) = \sum_{uv\in E(T)}\mathbb{H}(\varphi(u),\varphi(v)) - \sum_{v\in V(T)}(d_T(v)-1)\mathbb{H}(\varphi(v)).\]
\end{lem}

\begin{proof}
Let $k=e(T)$, let $v_0$ be an arbitrary ``root'' vertex of $T$ and label the other vertices of $T$ by $v_1,\dots,v_k$, in breadth-first-search order from $v_0$, so that the neighbours of $v_0$ in $T$ come first, followed by the vertices at distance two from $v_0$ in $T$, and so on. For each $1\leq i\leq k$, let $p(i)$ be the index $j<i$ such that $v_j$ is the neighbour of $v_i$ on the unique path from $v_i$ to $v_0$; we call $v_{p(i)}$ the \emph{parent} of $v_i$. By the Chain Rule (\cref{lem:ChainRule}), we have
\begin{align*}
\mathbb{H}(\varphi)&=\mathbb{H}(\varphi(v_0),\dots,\varphi(v_k))\\
&=\mathbb{H}(\varphi(v_0)) + \sum_{i=1}^k\mathbb{H}(\varphi(v_i)\mid \varphi(v_0),\dots,\varphi(v_{i-1})).
\end{align*}
Now, since $\varphi$ is a uniformly random homomorphism from $\vec{T}$ to $\vec{G}$, the random variables $\varphi(v_i)$ and $(\varphi(v_j): j<i, j\neq p(i))$ are conditionally independent given $\varphi(v_{p(i)})$. Indeed, once the image of $v_{p(i)}$ is fixed, the choices in the component of 
$T-v_iv_{p(i)}$ containing $v_i$ and the choices in the rest of the tree are 
constrained only through the image of $v_{p(i)}$. Thus, by \cref{lem:decon}, we have
\begin{align*}
    \mathbb{H}(\varphi) &= \mathbb{H}(\varphi(v_0)) + \sum_{i=1}^k\mathbb{H}(\varphi(v_i)\mid \varphi(v_{p(i)}))\\
    &= \mathbb{H}(\varphi(v_0)) + \sum_{i=1}^k\left(\mathbb{H}(\varphi(v_i), \varphi(v_{p(i)})) -\mathbb{H}(\varphi(v_{p(i)}))\right),
\end{align*}
where we have again applied the Chain Rule. For each edge $uv$ of $T$, either $\mathbb{H}(\varphi(u),\varphi(v))$ or the equivalent term $\mathbb{H}(\varphi(v),\varphi(u))$ appears in the above expression exactly once. For $j\geq1$, the above expression contains $\mathbb{H}(\varphi(v_j))$ with a coefficient of $-1$ exactly once for each index $i$ such that $p(i)=j$, of which there are $d_T(v_j)-1$. Also, $\mathbb{H}(\varphi(v_0))$ appears with a coefficient of $-1$ exactly $d_T(v_0)$ times and with a coefficient of $1$ exactly once. This completes the proof.
\end{proof}

Given a homomorphism $f$ from an oriented graph $\vec{H}$ to an oriented graph $\vec{G}$ and an arc ${(u,v) \in A(\vec{G})}$, let $f^{-1}((u,v))$ be the set of all arcs $(x,y)$ of $\vec{H}$ such that $(f(x),f(y))=(u,v)$, and let $f^{-1}(v)$ be the set of all vertices $x$ of $\vec{H}$ such that $f(x) = v$. We define $f^{-1}(uv)$ and $f^{-1}(v)$ similarly for a graph homomorphism $f$, edge $uv$ and vertex $v$. The next lemma is key to the proof of \cref{th:antiSidTreesStronger}. Similar ideas are used in, e.g.,~\cite[Lemma~2.6]{BehagueCrudeleNoelSimbaqueba25}.

\begin{lem}
\label{lem:coveringLemma}
Let $\vec{T}$ and $\vec{T}_1,\vec{T}_2,\dots,\vec{T}_m$ be (not necessarily distinct) oriented trees and let $f_i$ be a homomorphism from $\vec{T}_i$ to $\vec{T}$ for each $1\leq i\leq m$. If 
\begin{equation}\label{eq:edgeCover}\sum_{i=1}^m|f_i^{-1}((u,v))| = m,\qquad \forall (u,v)\in A(\vec{T})\end{equation}
and
\begin{equation}\label{eq:vertexCover}\sum_{i=1}^m|f_i^{-1}(v)| = m,\qquad \forall v\in V(\vec{T}),\end{equation}
then every oriented graph $\vec{G}$ satisfies
\[\hom(\vec{T},\vec{G})\leq \prod_{i=1}^m\hom(\vec{T}_i,\vec{G})^{1/m}.\]
\end{lem}

\begin{proof}
If $\hom(\vec{T},\vec{G})=0$, then the bound holds trivially. So, suppose that $\hom(\vec{T},\vec{G})>0$, and let $\varphi$ be a uniformly random homomorphism from $\vec{T}$ to $\vec{G}$. 

Let $\vec{H}$ be the disjoint union of $\vec{T}_1,\dots,\vec{T}_m$. Note that $\hom(\vec{H},\vec{G})=\prod_{i=1}^m\hom(\vec{T}_i,\vec{G})$. Our goal is to construct a distribution on the homomorphisms from $\vec{H}$ to $\vec{G}$ such that, if $\psi$ is chosen randomly with respect to this distribution, then $\mathbb{H}(\psi)=m\mathbb{H}(\varphi)$. Given that such a distribution exists, the result will follow since \cref{lem:uniform} implies that
\[\mathbb{H}(\psi)\leq \log_2(\hom(\vec{H},\vec{G})) =\log_2\left(\prod_{i=1}^m\hom(\vec{T}_i,\vec{G})\right)\]
and 
\[\mathbb{H}(\varphi)=\log_2(\hom(\vec{T},\vec{G})),\]
and so
\[m\log_2(\hom(\vec{T},\vec{G}))=m\mathbb{H}(\varphi)=\mathbb{H}(\psi)\leq\log_2\left(\prod_{i=1}^m\hom(\vec{T}_i,\vec{G})\right),\]
as desired.

Hence, it suffices to construct the desired distribution. To do so, we construct a random homomorphism $\psi$ from $\vec{H}$ to $\vec{G}$ via the following procedure. For each $1\leq i\leq m$, let $k_i=e(T_i)$. Start with an arbitrary root vertex $v_0^i$ of $\vec{T}_i$. Let $v_1^i,\dots,v_{k_i}^i$ be the other vertices of $\vec{T}_i$ listed in breadth-first-search order, that is, by distance from $v_0^i$. For each $1\leq j\leq k_i$, let $p_i(j)$ be the unique index such that $p_i(j)<j$ and one of $(v_j^i,v_{p_i(j)})$ and $(v_{p_i(j)},v_j^i)$ is an arc of $\vec{T}_i$. Let $\psi(v_0^i)$ be sampled according to the distribution of $\varphi(f_i(v_0^i))$. Then, for each $1\leq j\leq k_i$, sample $\psi(v_j^i)$ according to the distribution of $\varphi(f_i(v_j^i))$ conditioned on the event that $\varphi(f_i(v_{p_i(j)}^i))=\psi(v_{p_i(j)}^i)$, where this choice is conditionally independent of all previously chosen values of $\psi$ given $\psi(v_{p_i(j)}^i)$ (and choices for different trees $T_i$ are made independently). Applying the Chain Rule and conditional independence, we have
\begin{align*}
\mathbb{H}(\psi) &= \sum_{i=1}^m\left(\mathbb{H}(\psi(v_0^i)) + \sum_{j=1}^{k_i}\mathbb{H}(\psi(v_j^i)\mid\psi(v_{p_i(j)}^i))\right)\\
& = \sum_{i=1}^m\left(\mathbb{H}(\psi(v_0^i)) + \sum_{j=1}^{k_i}\left(\mathbb{H}(\psi(v_j^i),\psi(v_{p_i(j)}^i)) - \mathbb{H}(\psi(v_{p_i(j)}^i))\right)\right)\\
& = \sum_{i=1}^m\left(\sum_{xy\in E(T_i)}\mathbb{H}(\psi(x),\psi(y)) - \sum_{x\in V(T_i)}(d_{T_i}(x)-1)\mathbb{H}(\psi(x))\right).
\end{align*}
By construction, $(\psi(x),\psi(y))$ has the same distribution as $(\varphi(f_i(x)),\varphi(f_i(y)))$ for all ${xy\in E(T_i)}$, and $\psi(x)$ has the same distribution as $\varphi(f_i(x))$ for all $x\in V(T_i)$. So, the above expression can be rewritten as
\begin{align*}
\mathbb{H}(\psi) &= \sum_{i=1}^m\left(\sum_{xy\in E(T_i)}\mathbb{H}(\varphi(f_i(x)),\varphi(f_i(y))) - \sum_{x\in V(T_i)}(d_{T_i}(x)-1)\mathbb{H}(\varphi(f_i(x)))\right)\\
& = \sum_{i=1}^m \left(\sum_{uv\in E(T)}|f_i^{-1}(uv)|\cdot\mathbb{H}(\varphi(u),\varphi(v)) - \sum_{v\in V(T)}\sum_{x\in f_i^{-1}(v)}(d_{T_i}(x)-1)\mathbb{H}(\varphi(v))\right)\\
& = m\sum_{uv\in E(T)}\mathbb{H}(\varphi(u),\varphi(v)) - \sum_{v\in V(T)}\sum_{i=1}^m\sum_{x\in f_i^{-1}(v)}(d_{T_i}(x)-1)\mathbb{H}(\varphi(v))\\
& = m\sum_{uv\in E(T)}\mathbb{H}(\varphi(u),\varphi(v)) - \sum_{v\in V(T)}\sum_{i=1}^m\sum_{x\in f_i^{-1}(v)}d_{T_i}(x)\mathbb{H}(\varphi(v)) + m\sum_{v\in V(T)}\mathbb{H}(\varphi(v)),
\end{align*}
where the penultimate line follows from \eqref{eq:edgeCover} and the final line follows from \eqref{eq:vertexCover}. We now analyse the triple summation in the above expression, which counts, with multiplicity, incidences in the preimages of edges incident to $v$. By \eqref{eq:edgeCover}, for each $v\in V(T)$, each edge of $T$ incident to $v$ is ``covered'' exactly $m$ times by the homomorphisms $f_1,\dots,f_m$. Thus, since $f_1,\dots,f_m$ are homomorphisms, the sum of the degrees of vertices of $T_1,\dots,T_m$ that are mapped to $v$ by $f_1,\dots,f_m$ must be equal to $m\cdot d_T(v)$. Therefore,
\begin{align*}
    \mathbb{H}(\psi)&=m\sum_{uv\in E(T)}\mathbb{H}(\varphi(u),\varphi(v)) - \sum_{v\in V(T)}m\cdot d_{T}(v)\mathbb{H}(\varphi(v)) + m\sum_{v\in V(T)}\mathbb{H}(\varphi(v))\\
    &=m\mathbb{H}(\varphi)
\end{align*}
where the final line follows from \cref{lem:uniformHom}. This completes the proof. 
\end{proof}

We are now ready to prove \cref{th:antiSidTreesStronger}.

\begin{proof}[Proof of \cref{th:antiSidTreesStronger}]
Let $\vec{T}$ be a perfectly balanced oriented tree, constructed from an integer $k\geq0$, positive integers $d_1,\dots,d_k$ and the sets of vertices $S_0,S_1,\dots,S_k$. Let 
$\vec{T}_1,\dots,\vec{T}_{2^{e(T)}}$ enumerate all orientations of $T$, not 
identifying isomorphic copies.

We say that a homomorphism $f$ from $T$ to itself is \emph{level-preserving} if $f(S_i)\subseteq S_i$ for all $0\leq i\leq k$. Let $f_1,\dots,f_m$ enumerate each of the level-preserving homomorphisms from $T$ to $T$. In such a homomorphism, there are $|S_0|$ choices for the mapping of $S_0$ and there are $2d_i$ choices for the image of each vertex $v\in S_i$ for $1\leq i\leq k$. So, 
\[m=|S_0|\cdot \prod_{i=1}^k(2d_i)^{|S_i|}.\]
Observe that for each $j$ such that $1\leq j\leq m$, the homomorphism $f_j$ uniquely determines the orientation $\vec{T}^*$ of $T$ such that $f_j$ is a homomorphism from $\vec{T}^*$ to $\vec{T}$. Hence, for each $j$, there is a unique index $q(j)$ such that $f_j$ is a homomorphism from $\vec{T}_{q(j)}$ to $\vec{T}$. Also, for any $1\leq q\leq 2^{e(T)}$, the number of values of $j$ such that $q(j)=q$ is the same (in other words, the number of homomorphisms in $f_1,\dots,f_m$ that uniquely determine any given $\vec{T}_q$ is the same). Indeed, for any fixed orientation $\vec{T}_q$, in constructing such a homomorphism $f_j$, the mapping of $S_0$ is forced and there are $d_i$ choices for the image of each vertex $v\in S_i$ for each $1\leq i\leq k$. So, for each $q$, the number of values of $j$ such that $q(j)=q$ is 
\[\prod_{i=1}^kd_i^{|S_i|} = \frac{m}{2^{e(T)}}.\]
Since we have taken every level-preserving homomorphism of $T$, it is easy to see that $f_1,\dots,f_m$ satisfy the conditions of \cref{lem:coveringLemma}. Hence, applying \cref{lem:coveringLemma} with the homomorphisms $f_1,\dots,f_m$, and the oriented tree $\vec{T}_{q(i)}$ for each homomorphism $f_i$, we have
\begin{align*}
    \hom(\vec{T},\vec{G}) &\leq \prod_{i=1}^m \left(\hom(\vec{T}_{q(i)},\vec{G})\right)^{\frac{1}{m}}\\
    &= \prod_{i=1}^{2^{e(T)}} \left(\hom(\vec{T}_{i},\vec{G})^{\frac{m}{2^{e(T)}}}\right)^{\frac{1}{m}}\\
    &= \prod_{i=1}^{2^{e(T)}} \left(\hom(\vec{T}_{i},\vec{G})\right)^{\frac{1}{2^{e(T)}}},
\end{align*}
as desired.
\end{proof}

Finally, we show that \cref{th:antiSidTrees} follows from \cref{th:antiSidTreesStronger}.

\begin{proof}[Proof of \cref{th:antiSidTrees}]
Let $\vec{T}$ be a perfectly balanced oriented tree and let $\vec{T}_1,\dots,\vec{T}_{2^{e(T)}}$ enumerate all orientations of $T$, 
not identifying isomorphic copies. Let $\vec{G}$ be any oriented graph. Then, by \cref{th:antiSidTreesStronger} and the AM-GM inequality,
\begin{align*}
\hom(\vec{T},\vec{G})&\leq \left(\prod_{i=1}^{2^{e(T)}}\hom(\vec{T}_i,\vec{G})\right)^{1/2^{e(T)}}\\
&\leq \sum_{i=1}^{2^{e(T)}}\frac{\hom(\vec{T}_i,\vec{G})}{2^{e(T)}}\\
&=\frac{1}{2^{e(T)}}\hom(T,G)
\end{align*}
which completes the proof. 
\end{proof}

\subsection{Non-bipartite relatively anti-Sidorenko graphs}
\label{ssec:non_bpt_antiSid}

We conclude this section by proving \cref{th:antiSidChrom}, which we restate here for convenience.

\antiSidChrom*

\begin{figure}[htb]
\centering

\begin{tikzpicture}[
    scale=0.8,
    >={Stealth[scale=0.65]},
    font=\small,
    attachment/.style={->, line width=0.7pt, opacity=0.45},
    internal/.style={->}
]

\def\attachtoone#1#2#3{%
    \ifthenelse{\equal{#3}{p}}{\draw[attachment] (#1) -- (#2);}{\draw[attachment] (#2) -- (#1);}
}

\begin{scope}[xshift=0cm,yshift=0cm, scale=0.9]
    \node at (0,-0.75) {(a) $\vec{H}_2$};

    \node[mycircle] (one2) at (0,0) {};
    \node[mycircle,label=above left:{$-$}] (minus2) at (-1.45,1.25) {};
    \node[mycircle,label=above right:{$+$}] (plus2) at (1.45,1.25) {};

    \draw[->] (minus2) -- (one2);
    \draw[->] (one2) -- (plus2);
\end{scope}

\begin{scope}[xshift=7.6cm,yshift=0cm, scale=0.78]
    \node at (0,-4.5) {(b) $\vec{H}_3$};

    \node[mycircle] (one3) at (0,0) {};

    \def\R{3}
    \def\d{0.7}

    %
    %
    \foreach[count=\i from 0] \su/\sone/\sv/\lab in {
        p/p/p/{+++},
        p/p/m/{++-},
        p/m/p/{+-+},
        p/m/m/{+--},
        m/p/p/{-++},
        m/p/m/{-+-},
        m/m/p/{--+},
        m/m/m/{---}
    }{
        \pgfmathsetmacro{\ang}{90 - 45*\i}

        \node[mycircle] (U\i) at ({\R*cos(\ang) + \d*sin(\ang)},
                                  {\R*sin(\ang) - \d*cos(\ang)}) {};
        \node[mycircle] (W\i) at ({\R*cos(\ang)},
                                  {\R*sin(\ang)}) {};
        \node[mycircle] (V\i) at ({\R*cos(\ang) - \d*sin(\ang)},
                                  {\R*sin(\ang) + \d*cos(\ang)}) {};

        \draw[internal] (U\i) -- (W\i);
        \draw[internal] (W\i) -- (V\i);

        \attachtoone{one3}{U\i}{\su}
        \attachtoone{one3}{W\i}{\sone}
        \attachtoone{one3}{V\i}{\sv}

        \node[font=\tiny] at ({(\R+0.8)*cos(\ang)},{(\R+0.8)*sin(\ang)})
            {$\scriptstyle(\lab)$};
    }

\end{scope}

\end{tikzpicture}

\caption{
An illustration of the recursive construction in the proof of \cref{th:antiSidChrom} in the cases (a) $k=2$ and (b) $k=3$. 
}
\label{fig:non_bip_antisid_examples}
\end{figure}
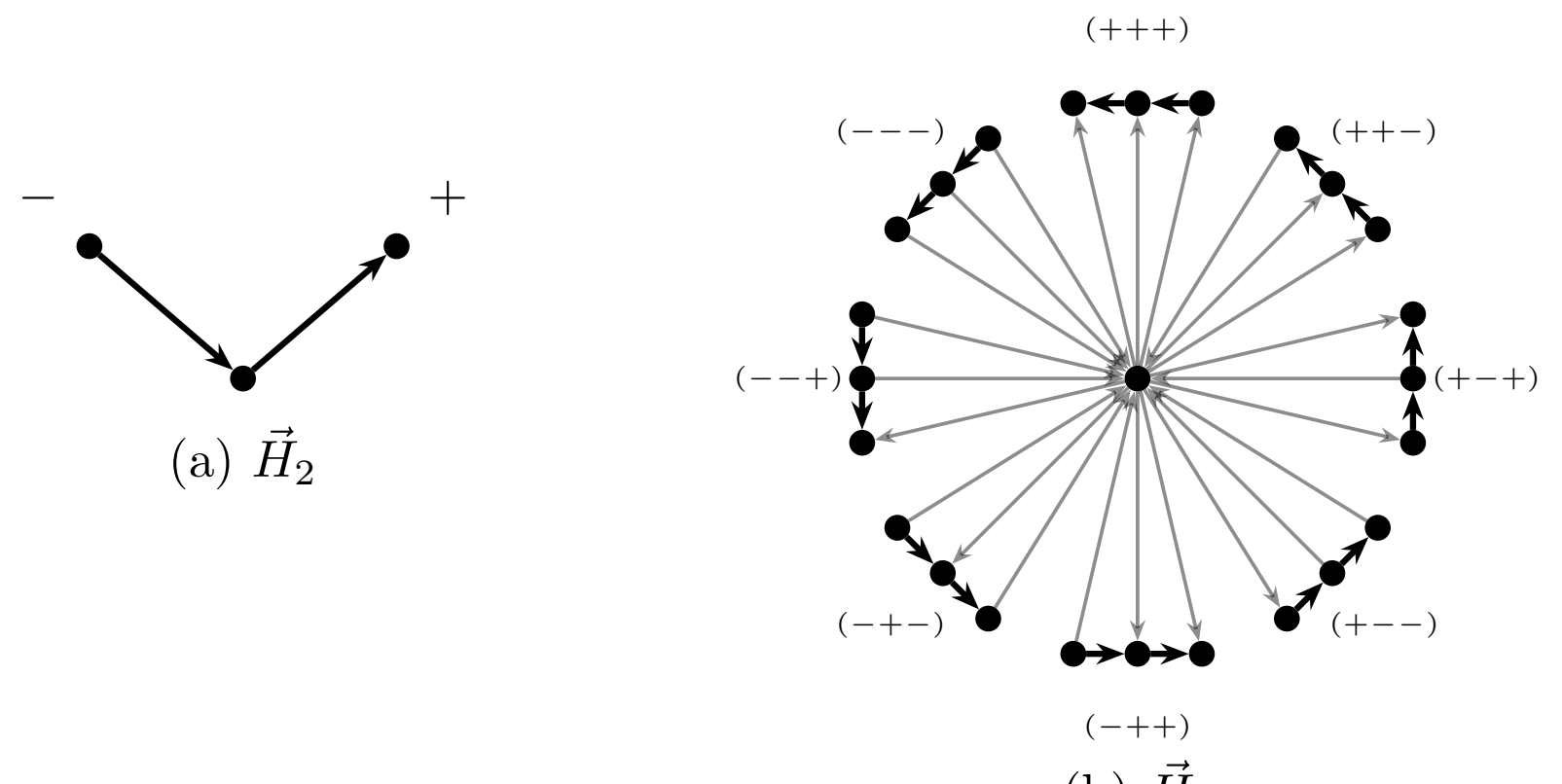

\begin{proof}
We show that there exists a relatively anti-Sidorenko oriented graph $\vec{H}_k$ such that $H_k$ has chromatic number $k$ for each $k\geq1$ by induction on $k$. For $k=1$, we simply let $\vec{H}_1$ consist of a single isolated vertex.

Now, suppose that $k\geq2$. By induction, we can assume that $\vec{H}_{k-1}$ is a relatively anti-Sidorenko oriented graph such that $\chi(H_{k-1})=k-1$. In constructing $\vec{H}_k$, we start with a single isolated vertex which we denote by $1$. Then, for each function $\sigma:V(\vec{H}_{k-1})\to\{-,+\}$, we add a disjoint copy $\vec{F}_\sigma$ of $\vec{H}_{k-1}$ where, for each $v\in V(\vec{H}_{k-1})$, we let $v_{\sigma}$ denote the copy of $v$ in $\vec{F}_{\sigma}$. For each $v\in V(\vec{H}_{k-1})$ and $\sigma:V(\vec{H}_{k-1})\to\{-,+\}$, add the arc $(1,v_{\sigma})$ if $\sigma(v)=+$ and the arc $(v_{\sigma},1)$ otherwise. This completes the construction of $\vec{H}_k$. See \cref{fig:non_bip_antisid_examples} for an illustration of $\vec{H}_k$ in the cases $k=2$ and $k=3$.

We start by showing that $\chi(H_k)=\chi(H_{k-1})+1=k$. Since vertex $1$ in $\vec{H}_k$ is adjacent to every vertex in each copy of $\vec{H}_{k-1}$, it follows that any proper colouring of $H_k$ assigns vertex $1$ to a colour that is distinct from the colours used on $H_{k-1}$, and so $\chi(H_k) \geq \chi(H_{k-1}) + 1$. By colouring each copy of $\vec{H}_{k-1}$ identically with $\chi(H_{k-1})$ colours, and introducing a new colour for vertex $1$, we get that $\chi(H_k) \leq \chi(H_{k-1}) + 1$ as well. 

So, all that remains is to prove that $\vec{H}_k$ is relatively anti-Sidorenko. Let $\vec{G}$ be an arbitrary oriented graph, and let $M = 2^{v(H_{k-1})}$. For $x\in V(G)$, let $N_G(x)$ denote the neighbourhood of $x$ in the underlying graph $G$. We start by counting the number of homomorphisms from $H_k$ to $G$. Such a homomorphism can be thought of as first selecting the image of vertex $1$, say $x$, and mapping each of the $M$ copies of $\vec{H}_{k-1}$ homomorphically to $G[N_G(x)]$. So, we have 
$$\hom(H_k,G) = \sum_{x\in V(G)}\hom(H_{k-1},G[N_G(x)])^M.$$
Since $\vec{H}_{k-1}$ is relatively anti-Sidorenko, $\hom(H_{k-1},G[N_G(x)])\geq 2^{e(H_{k-1})}\hom(\vec{H}_{k-1},\vec{G}[N_G(x)])$ for all $x\in V(G)$. Thus, 
\begin{equation}
\label{eq:applyAntiSid}
\hom(H_k,G)\geq 2^{M e(H_{k-1})}
   \sum_{x\in V(G)}
   \hom(\vec{H}_{k-1},\vec{G}[N_G(x)])^M.
\end{equation}
Now, for each pattern $\sigma : V(\vec{H}_{k-1}) \to \{+,-\}$ and $x\in V(G)$, define $h_\sigma(x)$ to be the number of homomorphisms $\phi : \vec{H}_{k-1} \to \vec{G}[N_G(x)]$ such that $(x,\phi(v))$ is an arc of $\vec{G}$ for all $v\in V(\vec{H}_{k-1})$ such that $\sigma(v)=+$ and $(\phi(v),x)$ is an arc of $\vec{G}$ for all $v\in V(\vec{H}_{k-1})$ such that $\sigma(v)=-$. Then, clearly, 
\[\hom(\vec{H}_{k-1},\vec{G}[N_G(x)]) = \sum_{\sigma}h_\sigma(x)\geq M\left(\prod_\sigma h_\sigma(x)\right)^{1/M}\]
where the sum and product are over all $\sigma : V(\vec{H}_{k-1}) \to \{+,-\}$ and the inequality follows from AM-GM. Substituting this into \eqref{eq:applyAntiSid} yields
\begin{align*}
\hom(H_k,G)&\geq 2^{M e(H_{k-1})}
   \sum_{x\in V(G)}
   M^M\left(\prod_\sigma h_\sigma(x)\right)\\
   &= 2^{e(H_{k})}
   \sum_{x\in V(G)}
   \left(\prod_\sigma h_\sigma(x)\right).
\end{align*}
Thus, to complete the proof, it suffices to argue that $\hom(\vec{H}_k,\vec{G})=\sum_{x\in V(G)}\left(\prod_\sigma h_\sigma(x)\right)$. This is true because, after selecting the image of vertex $1$, say $x$, the copy $\vec{F}_\sigma$ of $\vec{H}_{k-1}$ can be mapped into $\vec{G}$ in exactly $h_\sigma(x)$ ways, and these choices are independent of one another. Thus, the proof is complete. 
\end{proof}

\begin{rem}
\label{rem:antiSidConstruction}
The construction used in the proof of \cref{th:antiSidChrom} generalizes as follows. If $\vec{F}$ is a relatively anti-Sidorenko oriented graph and $\vec{H}$ is obtained from a vertex $1$ and $2^{v(\vec{F})}$ disjoint copies of $\vec{F}$ which are joined to vertex $1$ with arcs oriented in all of the possible ways, then $\vec{H}$ is relatively anti-Sidorenko as well. A similar construction in which one starts with two vertices $1$ and $2$ with an arc between them and adds $2^{2v(\vec{F})}$ copies of $\vec{F}$ joined to $\{1,2\}$ with arcs oriented in all of the possible ways also produces a relatively anti-Sidorenko oriented graph. 
\end{rem}

\section{Ruling Out Relatively Anti-Sidorenko Orientations}
\label{sec:notAntiSid}

We begin with a proof of \cref{prop:star}, as it does not require any sophisticated ideas. 
\antiSidStar*
\begin{proof}
First, we observe that an oriented star in which the central vertex has in-degree equal to its out-degree is a perfectly balanced oriented tree. Thus, all such oriented stars are relatively anti-Sidorenko by \cref{th:antiSidTrees}.

Therefore, let $\vec{S}_{a,b}$ be an oriented star with one central vertex and $a+b$ leaves, of which $a$ are sinks and $b$ are sources, with $a\neq b$ and $a+b\geq2$. We aim to show that $\vec{S}_{a,b}$ is not relatively anti-Sidorenko by counting homomorphisms into the directed host-star $\vec{S}_{Na,Nb}$, where $N$ is a large integer. We have $\hom(\vec{S}_{a,b},\vec{S}_{Na,Nb})=(Na)^a(Nb)^b$,
where we regard $0^0$ as being equal to $1$. Furthermore, letting $k=a+b$, we have
\[\hom(S_{a,b},S_{Na,Nb})=(N(a+b))^{a+b} + N(a + b) = N^kk^k + Nk.\]
Hence,
\begin{align*}
    \frac{\hom(\vec{S}_{a,b},\vec{S}_{Na,Nb})}{\hom(S_{a,b},S_{Na,Nb})}&=\frac{N^{k}a^ab^b}{N^{k}k^k + Nk} = \frac{a^ab^b}{k^k + N^{1-k}k}\\
    \implies \lim_{N\to\infty} \frac{\hom(\vec{S}_{a,b},\vec{S}_{Na,Nb})}{\hom(S_{a,b},S_{Na,Nb})} &=\frac{a^ab^b}{k^k}
\end{align*}
Note that the value of $a^ab^b$ is minimised when $a = b = k/2$, and strictly greater than this value for all other assignments of $a$ and $b$, subject to $a + b = k$. Hence, since $a \neq b$, we have that 
\begin{equation*}
    \frac{a^ab^b}{k^k} > \frac{(k/2)^k}{k^k} = \frac{1}{2^k},
\end{equation*}
and so for sufficiently large $N$, the host $\vec{S}_{Na,Nb}$ violates the relative anti-Sidorenko inequality, completing the proof.
\end{proof}

In the proof of \cref{prop:star}, we applied a strategy of finding many homomorphic copies of an oriented graph $\vec{H}$ in an oriented host graph $\vec{G}$ by taking $G$ to be a certain blow-up of $H$. Our approach to proving \cref{th:path} generalises this idea further. In particular, we will take each vertex and edge of $H$ to be weighted, and take a blow-up of $H$ where the size of each independent set corresponding to a given vertex of $H$ is determined by the weight of that vertex, and the density of the edges between two such vertex classes is determined by the weight of the corresponding edge in $H$. This approach is inspired by~\cite[Lemma~2.12]{BehagueCrudeleNoelSimbaqueba25} and is encapsulated by the following lemma. Given a graph $H$, a weight function $w:V(H)\cup E(H)\to[0,\infty)$ and a homomorphism $f$ from $H$ to itself, define 
\[w(f):=\sum_{v\in V(H)}w(f(v)) + \sum_{uv\in E(H)}w(f(u)f(v)).\]

\begin{lem}
\label{lem:dualLP}
Let $H$ be a graph. If there exists a weight function $w:V(H)\cup E(H)\to[0,\infty)$ such that every homomorphism $f$ from $H$ to itself satisfies
\begin{equation}\label{eq:w(f)}
w(f) \geq \sum_{v\in V(H)}w(v) + \sum_{e\in E(H)}w(e)
\end{equation}
and the number of homomorphisms $f$ from $H$ to itself which satisfy \eqref{eq:w(f)} with equality is less than $2^{e(H)}$, then $H$ does not admit a relatively anti-Sidorenko orientation. 
\end{lem}

\begin{proof}
We first show that we can assume $w$ is integer-valued. Given a choice of $w$ that satisfies the hypotheses of the lemma, we can let $\varepsilon$ be a small positive rational such that, for every homomorphism $f$ for which \eqref{eq:w(f)} is not tight, we have
\begin{equation}
\label{eq:w(f)2}
w(f) \geq \sum_{v\in V(H)}w(v) + \sum_{e\in E(H)}w(e) + \varepsilon.
\end{equation}
Thus, we can view $w$ as a feasible point for a system of inequalities that consists of \eqref{eq:w(f)} for every homomorphism $f$ from $H$ to $H$, \eqref{eq:w(f)2} for every such homomorphism for which \eqref{eq:w(f)} is not tight, along with non-negativity constraints on all vertex and edge weights. This is a system of inequalities in $|V(H)|+|E(H)|$ variables with rational coefficients and therefore it has a rational solution. Thus, we may assume that $w$ itself is rational. Multiplying by a common denominator, we can therefore assume that $w$ is integer-valued. Let 
\[m:=\sum_{v\in V(H)}w(v) + \sum_{e\in E(H)}w(e).\]

Now, fix an integer $k$ sufficiently large (with respect to $H$ and $w$), and choose $n$ large with respect to $H,w$ and $k$. Set $N=n^k$. We construct a graph $G$ as follows. For each vertex $v$ of $H$, create a set $S_v$ of $Nn^{-w(v)}$ vertices. For each $uv\in E(H)$ and each $x\in S_u$ and $y\in S_v$, include the edge $xy$ in $G$ with probability $n^{-w(uv)}$ independently of all other random choices. 

Note that there is a natural homomorphism $\pi : G \to H$ obtained by mapping each vertex of $S_v$ to $v$. We will calculate $\mathbb{E}[\hom(H,G)]$. For each homomorphism $f : H \to H$, let $X_f$ denote the set of homomorphisms $\phi : H \to G$ satisfying $\pi \circ \phi = f$. Then $\hom(H,G) = \sum_f |X_f|$, and so it suffices to calculate $\mathbb{E}[|X_f|]$ for each $f$. First, we observe that the expected number of injective elements of $X_f$ is at most 
\[\prod_{v\in V(H)}|S_{f(v)}|\cdot \prod_{uv\in E(H)}n^{-w(f(u)f(v))} =\prod_{v\in V(H)}Nn^{-w(f(v))}\cdot \prod_{uv\in E(H)}n^{-w(f(u)f(v))} = N^{v(H)}n^{-w(f)}.\]
Also, the number of non-injective elements of $X_f$ is (crudely) $O\left(N^{v(H)-1}\right)=O\left(N^{v(H)}n^{-k}\right)$. Recalling that $k$ is chosen large with respect to $H$ and $w$, we can ensure that $N^{v(H)-1}=o\left(N^{v(H)}n^{-w(f)}\right)$. Therefore, $\mathbb{E}[|X_f|]\leq(1+o(1))N^{v(H)}n^{-w(f)}$.

By the hypotheses of the lemma, every homomorphism $f$ from $H$ to itself satisfies ${w(f)\geq m}$ and at most $2^{e(H)}-1$ of them satisfy $w(f)=m$. Therefore,
\begin{equation}\label{eq:smallhomH}\mathbb{E}[\hom(H,G)] \leq \left(2^{e(H)}-1+o(1)\right)N^{v(H)}n^{-m},\end{equation}
where the $o(1)$ term tends to $0$ as $n$ tends to infinity. Now, consider any orientation $\vec{H}$ of $H$ and let $\vec{G}$ be the orientation of $G$ obtained by orienting the edges between $S_u$ and $S_v$ in the same direction as the orientation of $uv$ in $\vec{H}$. Then, by considering functions from $V(\vec{H})$ to $V(\vec{G})$ which map $v$ to $S_v$ for all $v\in V(\vec{H})$, we get that
\[\mathbb{E}[\hom(\vec{H},\vec{G})]\geq N^{v(H)}n^{-m}.\]
Combining this with \eqref{eq:smallhomH}, we see that $\hom(\vec{H},\vec{G})>(1/2^{e(H)})\hom(H,G)$ with positive probability and so $\vec{H}$ is not relatively anti-Sidorenko. This completes the proof.
\end{proof}

We close this section by applying \cref{lem:dualLP} to prove \cref{th:path}.

\notAntiSidPath*
\begin{proof}
Let $P_m$ denote the path with $m$ edges. First, consider an odd path $P_{2k+1}$ for some $k\geq1$. Let $w:V(P_{2k+1})\cup E(P_{2k+1})\to [0,\infty)$ be a weight function as follows.
\begin{align*}
    w(e) &= \begin{cases}
                0 & \text{ if $e$ is not incident to a leaf,}\\
                1 & \text{ if $e$ is incident to a leaf,}
            \end{cases}\\ 
    \intertext{and, }
    w(v) &= \begin{cases}
                0 & \text{ if $v$ is a leaf,}\\
                2 & \text{ if $v$ is not a leaf}.
            \end{cases}
\end{align*}

Then the sum of all weights is $2 + 2(2k)=4k+2$. Let $f$ be a homomorphism from $P_{2k+1}$ to itself. We bound $w(f)$ from below. Any vertex that is mapped to a non-leaf vertex contributes $2$ to $w(f)$, and for any vertex $v$ that is mapped to a leaf, the edges incident to $v$ each contribute $1$ to $w(f)$. Furthermore, since $f$ is a homomorphism and the length of the path is $2k+1\geq3$, no adjacent vertices can both be mapped to leaves and so the contributions from these edges are not being counted twice. Thus, letting $\ell$ be the number of leaves of $P_{2k+1}$ that $f$ maps to leaves of $P_{2k+1}$, we have 
\[w(f)\geq 2|V(P_{2k+1})| -\ell = 2(2k+2)-\ell = 4k+4-\ell.\]
Since $P_{2k+1}$ has only two leaves, we have $w(f)\geq 4k+2$ for all homomorphisms $f$. Moreover, $w(f)=4k+2$ only if $f$ maps the leaves of $P_{2k+1}$ to leaves; since the length of the path is odd, there are only two homomorphisms with this property. Since $2<2^{2k+1}$, it follows from \cref{lem:dualLP} that $P_{2k+1}$ has no relatively anti-Sidorenko orientation. 

Next, consider an even path $P_{2k}$ with $k\geq2$. Define $w:V(P_{2k})\cup E(P_{2k})\to [0,\infty)$ to be a slight modification of the weight function used for odd paths, as follows.
\begin{align*}
    w(e) &= \begin{cases}
                0 & \text{ if $e$ is not incident to a leaf}\\
                1 & \text{ if $e$ is incident to a leaf}
            \end{cases}\\ 
    \intertext{and}
    w(v) &= \begin{cases}
                0 & \text{ if $v$ is a leaf,}\\
                2 - \frac{1}{2k+2} & \text{ if $v$ is the central vertex,}\\
                2 & \text{ if $v$ is neither a leaf nor the central vertex}.
            \end{cases}
\end{align*}

Then the sum of all weights is $2+\left(2-\frac{1}{2k+2}\right) + 2(2k-2) = 4k-\frac{1}{2k+2}$. As in the case of paths of odd length, we analyse $w(f)$ for a homomorphism $f$ from $P_{2k}$ to itself by computing the contribution of each vertex. Every vertex contributes a weight of $2$ except for leaf vertices that are mapped to leaf vertices, which contribute $1$, and vertices that are mapped to the central vertex, which contribute $2-\frac{1}{2k+2}$. If we let $\ell$ be the number of leaf vertices that are mapped to leaf vertices and $c$ be the number of vertices mapped to the central vertex, then
\[w(f)=2|V(P_{2k})| - \ell - \frac{c}{2k+2} = 2(2k+1)-\ell-\frac{c}{2k+2} = 4k+2-\ell-\frac{c}{2k+2}.\]
If $f$ maps both leaves of $P_{2k}$ to leaves, then it can only map at most one vertex to the central vertex. So, if $\ell=2$, then $c\leq 1$. Hence, if $\ell=2$, then $w(f)\geq 4k-\frac{1}{2k+2}$ with equality if and only if exactly one vertex is mapped to the central vertex. Also, any homomorphism can only map at most $\frac{2k+2}{2}$ vertices to the central vertex, and therefore if $\ell\in\{0,1\}$, then $w(f)\geq 4k+1-\frac{1}{2}>4k-\frac{1}{2k+2}$. So, all homomorphisms from $P_{2k}$ to itself have weight at least $4k-\frac{1}{2k+2}$ and the only homomorphisms of weight equal to this are the four homomorphisms that map the leaves to the leaves and the central vertex to the central vertex. So, since $4<2^{2k}$, it again follows from \cref{lem:dualLP} that $P_{2k}$ has no relatively anti-Sidorenko orientation. 
\end{proof}

\section{Open Problems}
\label{sec:concl}
We conclude with several open problems inspired by the results of this paper. In \cref{sec:Sid} we show that the class of connected relatively Sidorenko graphs is precisely the class of one-way orientations of connected, reversible weakly norming graphs. However, a characterisation of the class of weakly norming graphs solely in terms of graph structure remains an open problem. Identifying a purely graph structural characterisation for the restricted class of reversible weakly norming graphs arising in this work would be an interesting step towards this problem, as well as providing new insights into the structure of our relatively Sidorenko graphs.

\begin{prob}
    Find a characterisation of the class of reversible weakly norming graphs entirely in terms of graph structure. Equivalently, find a characterisation of the class of connected relatively Sidorenko graphs in terms of graph structure.
\end{prob}

\cref{th:antiSidTrees} provides us with an infinite family of relatively anti-Sidorenko oriented trees, while~\cref{th:path} shows that not all trees admit a relatively anti-Sidorenko orientation. It would be interesting to have a characterisation of the complete class of relatively anti-Sidorenko trees.

\begin{prob}
\label{ques:antiSidTree}
Characterise the class of relatively anti-Sidorenko oriented trees.
\end{prob}

Using \cref{th:antiSidTrees}, one can obtain connected relatively anti-Sidorenko oriented graphs with $k$ arcs for any $k$ such that $k\not\equiv 3\bmod 4$. Moreover, the constructions described in \cref{ssec:non_bpt_antiSid} do not yield any connected relatively anti-Sidorenko oriented graphs with $k \equiv 3\bmod 4$. Interestingly, Nie and Wei recently obtained some AI-assisted results describing the class of relatively anti-Sidorenko trees in terms of a rooted variant of the same problem~\cite{NieWei2026personal}. These results imply that at least in the case of trees, this congruence condition is necessary. It remains unclear whether the congruence constraint $k \not\equiv 3 \bmod 4$ is required in general.

\begin{ques}
Does there exist a connected relatively anti-Sidorenko oriented graph $\vec{H}$ whose number of arcs is congruent to $3$ modulo $4$?
\end{ques}

Regarding the above questions, we believe that \cref{lem:dualLP} should be useful for ruling out the existence of relatively anti-Sidorenko orientations of a wide range of graphs. For example, using this lemma and linear programming software, we can verify that most trees on at most 8 vertices do not have such orientations. However, to rule out all oriented graphs apart from the known examples would require additional ideas.

\printbibliography

\end{document}